\documentclass[amstex,12pt]{amsart}
\usepackage{geometry}
\usepackage{mathtools,amssymb,amsthm,mathrsfs,color,lineno,paralist,graphicx,float}
\usepackage[linktocpage,backref=page,colorlinks,
linkcolor=magenta,
anchorcolor=magenta,
citecolor=magenta,
]{hyperref}
\usepackage{color}
\usepackage[T1]{fontenc}
\usepackage[utf8]{inputenc}
\usepackage{mathrsfs}
\usepackage{cleveref}
\usepackage{mathscinet}
\usepackage[normalem]{ulem}
\usepackage{bm}

\usepackage{calc}
\newcommand{\graph}{\text{graph}}

\newcommand{\diam}{\text{diam}}
\usepackage{etoolbox}

\newcounter{constK}

\newtheorem{thm}{Theorem}[section]
\theoremstyle{plain}
\theoremstyle{definition}
\newtheorem{proposition}[thm]{Proposition}

\newtheorem{lemma}[thm]{Lemma}
\newtheorem{definition}[thm]{Definition}
\newtheorem{claim}[thm]{Claim}

\newtheorem{remark}[thm]{Remark}

\numberwithin{equation}{section}

\definecolor{lgray}{gray}{0.9}

\def\leq{\leqslant}
\def\geq{\geqslant}
\def\tilde{\widetilde}
\def\hat{\widehat}

\allowdisplaybreaks
\renewcommand*{\backref}[1]{}
\renewcommand*{\backrefalt}[4]{\quad \tiny
  \ifcase #1 (\textbf{NOT CITED.})
  \or    (Cited on page~#2.)
  \else   (Cited on page~#2.)
  \fi}
\makeatletter

\title[Local entropy for non-compact systems and its application]{The Brin-Katok local entropy formula for non-compact random dynamical systems and an application}

\author[Z. Lian]{Zeng Lian}
\address[Z. Lian]
{	School of Mathematics,
Sichuan University,
Chengdu, Sichuan, 610016, China} \email{lianzeng@scu.edu.cn, zenglian@gmail.com}

\author[J. Zhang]{Jianhua Zhang}
\address[J. Zhang]{School of Mathematical Sciences, 
Chengdu University of Technology, Chengdu, Sichuan, 610059, China}
\email{leapforg@mail.ustc.edu.cn}

\date{\today}

\begin{document}

\begin{abstract}
   In this paper, we establish the classical Brin-Katok  local entropy formula for  $C^1$ random dynamical systems   on a  complete Riemannian manifold (possibly non-compact) under a mild integrability condition.   The proof depends on the constructions of  measurable partitions adapted to asymptotically varying radii and  finite measurable  partitions with null-boundaries.  As an application, under additional integrability assumptions, we obtain the Margulis-Ruelle inequality for $C^2$ random dynamical systems on complete non-compact Riemannian manifolds. Furthermore,  we also construct a counterexample showing that the integrability condition above cannot in general be omitted.   
        \end{abstract}
\maketitle

\tableofcontents

\medskip

\section{Introduction}
\label{24-4-13-2040}

Measure-theoretic entropy, originating from the pioneering work of Kolmogorov and Sinai \cite{MR0103254, MR103256}, serves as a fundamental concept for quantifying the complexity of measure-preserving dynamical systems. Since its inception, various formulations of entropy have been introduced, and their interrelationships have been systematically established through results such as the variational principle, Brin–Katok local entropy formula, and Katok's entropy formula \cite{MR175106, MR289746, MR730261, MR573822, MR274707}.

In this paper, we focus on the relation between local entropy and measure-theoretic entropy.  In 1983, Brin and Katok \cite{MR730261} established the following local entropy formula: let $X$ be a compact metric space and $T: X \to X$ be a continuous map. For any $T$-invariant probability measure $\mu$, we have for $\mu$-almost every $x \in X$:
  \begin{align}
 \label{26-3-14-1628}
 h_{\mu}(T,x):=\lim_{\epsilon\to 0}\liminf_{n\to\infty}\frac{-\log\mu(B(x,n,\epsilon))}{n}=\lim_{\epsilon\to 0}\limsup_{n\to\infty}\frac{-\log\mu(B(x,n,\epsilon))}{n},
 \end{align}
 where $B(x,n,\epsilon)=\{y\in X: d(T^ix, T^iy)<\epsilon, i=0,1,\cdots,n-1\}$ is the Bowen ball. The measure-theoretic entropy is recovered via $h_{\mu}(T) = \int_X h_{\mu}(T,x)\mathrm{d}\mu(x)$. This formula has since been extended to various settings, including conditional entropy and random dynamical systems \cite{MR3485864, MR3834712, MR2511539, MR3956171}.

Non-compact phase spaces arise naturally in smooth and random dynamical systems, for instance in systems on Euclidean spaces, open manifolds, and state spaces of stochastic differential equations. Entropy theory for non-compact systems has been studied in various settings; see, for example, \cite{MR1348316,MR2718907,MR3809061,MR3982746} and the references therein. A basic difference from the compact case is that invariant probability measures are not automatic. In the compact deterministic setting, the Krylov--Bogolyubov argument yields invariant probability measures for continuous maps, while in non-compact spaces empirical measures may fail to be tight and mass may escape to infinity. Existence results in non-compact settings usually require additional tightness, dissipativity, or recurrence assumptions ( see \cite{A,DaPratoZabczyk1996,MeynTweedie1993}).

In this paper,  we don't  address the existence problem for invariant measures. Instead, we assume that an invariant probability measure is given, and ask whether its entropy can still be recovered from local orbit statistics.  The main difficulty is that the geometry at infinity may affect Bowen balls and measurable partitions. The counterexample in \Cref{26-5-26-0943} shows that the Brin--Katok formula may fail  without additional assumptions. To control this behavior, we impose the integrability condition \eqref{2026-2-15-1654} on the local expansion of the system. Under this condition, we prove that the classical Brin--Katok local entropy formula remains valid for random dynamical systems on complete Riemannian manifolds; see \Cref{26-2-13-1805}. The result appears to be new already in the deterministic non-compact case, namely when the driving space consists of a single point.

The proof of \Cref{26-2-13-1805} has two main parts. To show that the upper local entropy is bounded above by the measure-theoretic entropy, we introduce a variant of Brin--Katok local entropy associated with an asymptotic sequence of radii. A main technical device is to allow the radii of Bowen balls to vary along return times, which makes it possible to construct measurable partitions of finite conditional entropy adapted to the non-compact geometry. Combining this construction with a conditional Shannon--McMillan--Breiman theorem for the skew product gives the desired upper estimate. For the reverse inequality, we use the regularity of the measure to construct finite measurable partitions with null boundaries and adapt Katok's argument \cite{MR573822} to the random non-compact setting.


As an application of the Brin--Katok local entropy formula, we prove a Margulis--Ruelle inequality for $C^2$ random dynamical systems on complete non-compact Riemannian manifolds. The  Margulis--Ruelle inequality bounds the measure-theoretic entropy from above by the integral of the sum of positive Lyapunov exponents, counted with multiplicities, and is well understood in compact smooth settings \cite{MR1314494,MR2005856,MR516310}. In the non-compact case, however, the inequality may fail without additional control at infinity \cite{MR3612326}. The Margulis-Ruelle inequality for deterministic non-compact systems has been obtained under suitable geometric or integrability assumptions \cite{MR4417345,MR4973346,MR3809061}. Under the assumptions \hyperref[23-4-29-1706]{(H1)}--\hyperref[24-9-26-1718]{(H2)}, we obtain in \Cref{24-4-12-1500} a random version of the inequality. The proof estimates the volume of neutralized Bowen balls in Lyapunov charts, using ideas related to \cite{MR4874980,MR4756121}.

\subsection{Brin-Katok local entropy formula}
In this subsection, we state the precise form of the Brin--Katok local entropy formula in our setting. Throughout the  paper, we  assume that $M$ is a complete  Riemannian manifold and $(\Omega,\mathscr{F},\mathbb{P},\theta)$ is an invertible  metric system\footnote{Throughout this paper, a metric system means that it is a measure-preserving dynamical system on a Polish probability space.}. We first recall the basic definitions and notations needed to state the main results. For further details, the reader is referred to \cite{Lian_10,MR2823209,MR1369243,MR3595937,MR3663618}.
\begin{definition}[Random dynamical systems]
    A  random dynamical system (abbreviated as RDS)  on $M$ over the metric system $(\Omega, \mathscr{F}, \mathbb P, \theta)$ is a measurable map 
    \[\varphi: \mathbb{Z}_+\times\Omega\times M\to M\] 
    satisfying $\varphi(0,\omega)=\mathrm{id}_{M}$ for all $\omega\in \Omega$ and the following cocycle property: 
    \begin{align*}
        \varphi(m+n,\omega) = \varphi(n,\theta^m \omega)\circ\varphi(m,\omega), \quad \text{for all } m,n\in \mathbb{Z}_+, \omega \in \Omega. 
    \end{align*}
   Given an integer $k\geq 0$,   we say the  RDS $\varphi$ is $C^k$,  if for every $(n,\omega)\in \mathbb{Z}_+\times \Omega$, the mapping 
    \[\varphi(n,\omega)\in C^k(M,M).\]
In particular, a $C^0$ RDS is called a continuous RDS.
  \end{definition}

A typical example of a   RDS related to SPDEs or SDEs is the one generated by the time-one maps of solutions (for example, see \cite{MR3443633}). In fact, the RDS  $\varphi: \mathbb{Z}_{+} \times \Omega \times M \rightarrow M$  can be written as follows
\begin{align*}
    (n,\omega, x)\mapsto \varphi^n_{\omega}(x):=\varphi(n,\omega)x=
    \begin{cases}
        \varphi_{\theta^{n-1}\omega}\circ\cdots \circ \varphi_{\omega}(x)&\text{ if } n\geq 1,
        \\ x &\text{ if } n=0, 
    \end{cases}
\end{align*}
by letting  $\varphi_{\omega}:=\varphi(1,\omega)$.  The RDS $\varphi$ generates  a measurable \emph{skew-product} transformation on $\Omega\times M$:
    \begin{align}
    \label{23-1-3-1520}
    \varPhi: \Omega\times M\to \Omega\times M\quad\text{with}\quad (\omega, x)\mapsto (\theta \omega, \varphi_{\omega}(x)),
    \end{align}
 which defines a measurable dynamical system on $(\Omega\times M, \mathscr{F}\otimes\mathscr{B}(M))$.  

\begin{definition}[Invariant measures for RDS]
    A Borel probability measure $\mu$ on $\Omega \times M$ is called an \emph{invariant measure} (resp. \emph{ergodic measure}) of an RDS $\varphi$ if it satisfies the following:
    \begin{enumerate}[(a)]
        \item $\mu$ is invariant (resp. ergodic) under the skew-product transformation $\varPhi$ associated with $\varphi$ as defined in \eqref{23-1-3-1520}; 
        \item $(\pi_{\Omega})_*\mu = \mathbb{P}$, where $\pi_{\Omega}$ denotes the projection from $\Omega \times M$ to $\Omega$. 
    \end{enumerate}
    Such a pair $(\varphi, \mu)$ is called a (resp. ergodic) \emph{measure-preserving RDS}.
\end{definition}
Usually, the  Borel  probability measure $\mu$ on $\Omega\times M$ with $(\pi_{\Omega})_*\mu=\mathbb{P}$  is characterized by its sample measures $\{\mu_{\omega}\}_{\omega\in\Omega}$ on $M$  with the  following properties 
  \begin{itemize}
    \item for any bounded measurable $f:M\to\mathbb{R}$,  the mapping $\omega\mapsto \int_{M}f\mathrm{d}\mu_{\omega}$ is $\mathscr{F}$-measurable;
    \item for any bounded measurable $g:\Omega\times M\to\mathbb{R}$, one has that 
    $$\int_{\Omega\times M}g\,\mathrm{d}\mu=\int_{\Omega}\int_{M}g\,\mathrm{d}\mu_{\omega}\mathrm{d}\mathbb{P}(\omega).$$
    \end{itemize}
    It  is written symbolically $\mu=\int_{\Omega}\mu_{\omega}\mathrm{d}\mathbb{P}$, called the disintegration of $\mu$ with respect to the marginal $\mathbb{P}$. The reader can see the details in  \cite[Proposition 3.6]{MR1993844}.


\medskip

We now state the Brin--Katok local entropy formula.
\begin{thm}
\label{26-2-13-1805}
 Let  $(\varphi,\mu)$ be a $C^1$ measure-preserving RDS on $M$ over an  invertible  metric system $(\Omega,\mathscr{F},\mathbb{P},\theta)$. If it satisfies the following integrability condition  
    \begin{align} 
    \label{2026-2-15-1654}
        \int_{\Omega\times M}\Big(\log^+\sup_{y\in B(x,1)}|\mathrm{d}_{y}\varphi_{\omega}|\Big)\mathrm{d}\mu(\omega,x)<+\infty,
        \end{align}
 then  the Brin-Katok local entropy formula holds. Namely,  for $\mu$-a.s. $(\omega,x)\in \Omega\times M$ one has that 
    \begin{align}
          h_{\mu}(\varphi)(\omega,x):=\lim_{\epsilon\to0}\liminf_{n\to\infty}-\frac{\log\mu_{\omega}(B(\omega,x, n, \epsilon))}{n}=\lim_{\epsilon\to0}\limsup_{n\to\infty}-\frac{\log\mu_{\omega}(B(\omega,x, n,\epsilon))}{n},
    \end{align}
    where $B(\omega,x,n,\epsilon):=\{y\in M: d(\varphi^i_{\omega}x,\varphi^i_{\omega}y)<\epsilon\text{ for } i=0,1,\dots,n-1\}$,  and the measure-theoretic entropy of RDS $(\varphi,\mu)$ is equal to 
    $$h_{\mu}(\varphi)=\int_{\Omega\times M} h_{\mu}(\varphi)(\omega,x)\mathrm{d}\mu(\omega,x).$$ 
     \end{thm}
     
\begin{remark}
In fact, the assumption in \Cref{26-2-13-1805} can be weakened from $C^1$ to locally Lipschitz. Particularly, \eqref{2026-2-15-1654} can be replaced with 
$$ \int_{\Omega\times M}\Big(\log^+\sup_{y\neq z\in B(x,1)}\frac{d\big(\varphi_{\omega}(y), \varphi_{\omega}(z)\big)}{d(y, z)}\Big)\mathrm{d}\mu(\omega,x)<+\infty.$$
The stronger $C^1$ condition was adopted merely to simplify the presentation.  However, the assumption \eqref{2026-2-15-1654} cannot in general be omitted. See  \Cref{26-5-26-0943} for a counterexample.
\end{remark}

     \subsection{Margulis-Ruelle inequality}
  To obtain  Margulis-Ruelle inequality on non-compact Riemannian manifolds, we impose geometric control and integrability assumptions on the first and second derivatives.   First,  we  give some  geometric estimations of the exponential map.  The following lemma follows from standard Riemannian geometry, and the reader can refer to \Cref{26-7-26-1608} for a detail proof.
\begin{lemma}
\label{26-1-22-1459}
 Given a small positive constant $b\ll 1$,  there exists a measurable  \emph{regular radius}  $\rho_{b}: M\to (0,1)$ such that for any $x\in M$
\begin{itemize}
\item $\exp_x: \mathcal{B}_x(\rho_b(x))\to B(x,1)$ is a diffeomorphism onto its image, where   $\mathcal{B}_x(\rho_b(x))$ is the open ball in  $T_xM$ at the  origin with  radius $\rho_b(x)$;
\item  $|\mathrm{d}_{.}\exp_x|\in (1-b, 1+b)$ and $|\nabla_{.}(\mathrm{d}\exp_x)|\leq b $ on $\mathcal{B}_x(\rho_b(x))$.
\end{itemize}
\end{lemma}
Note that if $M=\mathbb{R}^n$, then we can take the regular radius $\rho_b\equiv b$.
Let $(\varphi,\mu)$ be a $C^2$ measure-preserving RDS on $M$ over $(\Omega,\mathscr{F},\mathbb{P},\theta)$.   We  impose the following  mild integrability conditions on $(\varphi,\mu)$
 \begin{align} 
  \tag{H1}  \label{23-4-29-1706}
        &\int_{\Omega\times M}\Big(\log^+\sup_{y\in B(x,1)}|\mathrm{d}_{y}\varphi_{\omega}|+\log^+\sup_{y\in B(x,1)}|\nabla_{y}\mathrm{d}\varphi_{\omega}|\Big)\mathrm{d}\mu(\omega,x)<+\infty,
        \\\tag{H2}\label{24-9-26-1718}&\int_{\Omega\times M}|\log\rho_{b}(x)|\mathrm{d}\mu(\omega,x)<+\infty,
    \end{align}
    where $\rho_b$ is the regular radius defined as \Cref{26-1-22-1459} for some positive constant $b>0$.   
    \begin{thm}
\label{24-4-12-1500}
    For a $C^2$   measure-preserving  RDS $(\varphi,\mu)$ on $M$ over an  invertible  metric system $(\Omega,\mathscr{F},\mathbb{P},\theta)$, if it satisfies  {\rm \hyperref[23-4-29-1706]{(H1)}-\hyperref[24-9-26-1718]{(H2)}},      then the Margulis-Ruelle  inequality  holds:
    \begin{align}
        \label{24-7-21-2338}
        h_{\mu}(\varphi)\leq\int_{\Omega\times M}\sum_{\lambda_i>0}m_i(\omega,x)\lambda_i(\omega,x)\mathrm{d}\mu(\omega,x),
    \end{align}
    where $\{(\lambda_i,m_i)\}$ is the  Lyapunov spectrum of $(\varphi,\mu)$. 
     \end{thm}

The paper is organized as follows. In \Cref{24-4-12-1442}, we recall measure-theoretic entropy and several forms of local entropy for random dynamical systems. In \Cref{26-3-15-1525}, we construct measurable partitions adapted to asymptotic sequences and prove that the upper Brin--Katok local entropy is bounded above by the measure-theoretic entropy. In \Cref{24-4-12-1445}, we construct finite partitions with null boundaries and prove the reverse inequality. In \Cref{26-3-15-1526}, we prove the Margulis--Ruelle inequality by estimating the volume of neutralized Bowen balls in Lyapunov charts. Finally, in \Cref{26-5-26-0943}, we give a counterexample showing that the integrability condition \eqref{2026-2-15-1654} cannot in general be omitted.


\section{Entropies for  random dynamical systems}
\label{24-4-12-1442}
In this section, we recall the definition and basic properties of measure-theoretic entropy for random dynamical systems. We then introduce several variants of Brin--Katok local entropy that will be used in the proof of the main results. Throughout this section, $M$ is a complete Riemannian manifold and $\varphi$ is a continuous RDS on $M$ over an invertible metric system $(\Omega,\mathscr{F},\mathbb{P},\theta)$.
\subsection{Measure-theoretic entropy for RDS}
In this subsection, we will successively  give the definition of the measure-theoretic  entropy of  a measure-preserving RDS $(\varphi,\mu)$. 
 For any finite or countable Borel measurable partition $\alpha$ of $\Omega\times M$,  for any $(\omega,x)\in\Omega\times M$,  denote $\alpha(\omega,x)$ as the unique element in $\alpha$ containing the point $(\omega,x)$, and denote   
 $$\alpha_{\omega}:=\left\{ A_{\omega}: A\in \alpha\right\},$$
as the Borel measurable partition on $M$ induced by $\alpha$, where $A_{\omega}:=\{y\in M: (\omega,y)\in A\}$ is fiber of $A$  at $\omega$.   For a Borel probability measure  $\nu$ on $M$ and a Borel measurable partition $\beta$ of $M$, we set $H_{\nu}(\beta):=\sum_{A\in\beta}- \nu(A)\log\nu(A)$.   Hence   the \emph{conditional static entropy} of the partition $\alpha$  is defined  as 
$$H_{\mu}(\alpha|\pi_{\Omega})=\int_{\Omega}H_{\mu_{\omega}}(\alpha_{\omega})\mathrm{d}\mathbb{P}(\omega),$$ 
where $\pi_{\Omega}$ denotes the projection from $\Omega \times M$ to $\Omega$. 
The \emph{entropy} of the partition  $\alpha$  with respect to $(\varphi,\mu)$ is defined as 
\begin{align}
h_{\mu}(\varphi,\alpha):=\lim_{n\to+\infty}\frac{1}{n}H_\mu\left(\bigvee_{i=0}^{n-1}\varPhi^{-i}\alpha\big|\pi_{\Omega}\right),
\end{align}
where $\varPhi$ is the associated skew-product transformation defined as    \eqref{23-1-3-1520}.
\begin{definition}[Measure-theoretic entropy of RDS]
\label{25-11-11-1122}
    The \emph{measure-theoretic  entropy} of  a measure-preserving RDS $(\varphi,\mu)$ on $M$ is defined as
    $$h_\mu(\varphi):=\sup_{\alpha}h_{\mu}(\varphi,\alpha),$$
   where  $\alpha$ runs over countable Borel measurable partitions of $\Omega\times M$ with finite conditional static entropy. 
\end{definition} 
\begin{remark}
In fact, the measure-theoretic entropy of RDS coincides with the conditional measure-theoretic entropy of factor  map 
$$\pi_{\Omega}: (\Omega\times M,\mathscr{F}\otimes\mathscr{B}(M),\mu, \varPhi)\to (\Omega,\mathscr{F},\mathbb{P},\theta).$$
For the proof, we refer the reader to \cite{bogenschutz1993equilibrium}.
\end{remark}

Finally, we will give  another equivalent form of measure-theoretic  entropy of RDS.  The reader can refer to the proof in \cite[Theorem 4.5 in Chapter 0]{MR1369243}.
\begin{proposition}
\label{26-2-22-2151}
Let $(\varphi,\mu)$ be a measure-preserving RDS on $M$ over a metric system $(\Omega,\mathscr{F},\mathbb{P},\theta)$.   If $\{\alpha_n\}_{n\in\mathbb{N}}$ is an increasing finite measurable partition of $M$ such that 
$$\pi_{M}^{-1}\Big(\bigvee_{n\in\mathbb{N}}\hat{\alpha}_n\Big)\vee\pi_{\Omega}^{-1}(\mathscr{F})=\mathscr{F}\otimes\mathscr{B}(M)\pmod\mu,$$
where $\pi_M$ is coordinate projection from $\Omega\times M$  to $M$,  and $\hat{\alpha}_n$ is  $\sigma$-algebra generated by the partition $\alpha_n$, then one has that 
$$h_{\mu}(\varphi)=\lim_{n\to\infty}h_{\mu}\big(\varphi, \pi^{-1}_{M}(\alpha_n)\big).$$
 \end{proposition}
\subsection{Local entropies for RDS}
\label{24-11-2-1624}
In this subsection, we introduce  two   modified version of the Brin-Katok local entropies of measure-preserving RDS $(\varphi,\mu)$ and establish some of their basic properties.  These  concepts  are inspired  by \cite{MR3233706,MR4756121}.

  Let $\mathbf r:\Omega\times M\to\mathbb R_+$ be a Borel measurable function, and let $A\subset \Omega\times M$ be a measurable set with positive $\mu$-measure. For $n\in\mathbb N$ and $(\omega,x)\in\Omega\times M$, define the modified Bowen ball

\begin{align}
\label{24-8-25-1035}
    B_{A}(\omega,x, n,\mathbf{r})=\Big\{y\in M:\,&\varphi^{i}_{\omega}y\in B\big(\varphi^{i}_{\omega}x, \mathbf{r}(\varPhi^{i}(\omega,x))\big)\cap A_{\theta^i\omega}
    \\\notag\quad&\text{ whenever  } i\in\{0,\dots,n-1\}\text{ and }\varPhi^{i}(\omega,x)\in A\Big\}.
\end{align}


\begin{definition}[Asymptotic sequence]\label{d.010901}
A family of Borel measurable functions $\mathfrak{r}:=\{\mathbf{r}_{m}: \Omega\times M\to\mathbb{R}_+\}_{m\in\mathbb{N}}$
is called an \emph{asymptotic sequence} on $\Omega\times M$ if it satisfies the following properties for $\mu$-a.s. $(\omega,x)\in \Omega\times M$:
        \begin{enumerate}[(i)]
        \item  $\log\mathbf{r}_{m}$ is $\mu$-integrable for each $m\in\mathbb{N}$;
            \item $\mathbf{r}_{m_1}(\omega,x)\leq \mathbf{r}_{m_2}(\omega,x)$  for any $m_1>m_2$;
           \item  $\lim_{m\to\infty}\mathbf{r}_{m}(\omega,x)=0.$
        \end{enumerate}
In this paper, we shall consider asymptotic sequences that depend on a given positive $\mu$-measurable subset $A\subset \Omega\times M$. Accordingly, we shall also refer to the pair $(\mathfrak{r}, A)$ as an asymptotic sequence.
\end{definition}
\begin{definition}[Local entropies]\label{24-9-6-1111}
    Given a general measure-preserving RDS $(\varphi,\mu)$ on $M$ over the metric system $(\Omega,\mathscr{F},\mathbb{P},\theta)$, its \emph{ upper  Brin-Katok local entropy} at $(\omega,x)\in\Omega\times M$ with respect to an asymptotic sequence $(\mathfrak{r},A)$  is defined as 
    \begin{align*}
        &\overline{h}^{(\mathfrak{r}, A)}_{\mu}(\varphi)(\omega,x)=\lim_{m\to\infty}\limsup_{n\to\infty}-\frac{1}{n}\log\mu_{\omega}(B_{A}(\omega,x, n, \mathbf{r}_{m})).
    \end{align*}
     If we take $\mathbf{r}_m\equiv1/m$ and $A=\Omega\times M$,  then  \emph{upper Brin-Katok local entropy} with respect to  $(\mathfrak{r},A)$ is the upper classical Brin-Katok entropy, which is denoted as 
   \begin{align*} 
    &\overline{h}_{\mu}(\varphi)(\omega,x)=\overline{h}^{(\mathfrak{r}, A)}_{\mu}(\omega,x)=\lim_{m\to\infty}\limsup_{n\to\infty}-\frac{1}{n}\log\mu_{\omega}(B(\omega,x,n, 1/m)).
    \end{align*}
  Denote the \emph{neutralized  Bowen ball} as 
  $$B(\omega,x,n, e^{-n\epsilon}):=\{y\in M: \varphi^i_{\omega}y\in B(\varphi^i_{\omega}x, e^{-n\epsilon})\text{ for }i=0,1\dots,n-1\}.$$
    Then   \emph{upper neutralized Brin-Katok local entropy} of $(\varphi,\mu)$ at $(\omega,x)\in\Omega\times M$   is defined as 
    \begin{align}\label{e.010601}
        \overline{h}^{N}_{\mu}(\varphi)(\omega,x):=\lim_{m\to\infty}\limsup_{n\to\infty}-\frac{1}{n}\log\mu_{\omega}(B(\omega,x, n, e^{-n/m})).
    \end{align}
    
\end{definition}

 \begin{remark}
 Note that   for any $p>0$ and $(\omega,x)\in\Omega\times M$ one has that 
     \begin{align*}
     \overline{h}^{(\mathfrak{r}, A)}_{\mu}(\varphi)(\omega,x)&=\lim_{m\to\infty}\limsup_{n\to\infty}-\frac{1}{n}\log\mu_{\omega}(B_{A}(\omega,x, n, p\mathbf{r}_{m})),
  \\\overline{h}^{N}_{\mu}(\varphi)(\omega,x)&=\lim_{m\to\infty}\limsup_{n\to\infty}-\frac{1}{n}\log\mu_{\omega}(B(\omega,x, n, pe^{-n/m})).
     \end{align*}
 \end{remark}

Now, we are going to  prove that the upper  Brin-Katok local entropy function $(\omega,x)\mapsto \overline{h}^{(\mathfrak{r},A)}(\varphi)(\omega,x)$ is Borel measurable by applying the following lemma. 
\begin{lemma}
\label{24-10-27-1252}
    Let 
    $$A: \Omega\times M\to 2^{M},\quad(\omega,x)\mapsto A(\omega,x)$$ be a set-valued map such that  
    $$\text{graph}(A):=\{(\omega,x,y)\in \Omega\times M\times M: y\in A(\omega,x)\}$$ 
 is a measurable subset of $\Omega\times M\times M$. Then the mapping $(\omega,x)\mapsto\mu_{\omega}(A(\omega,x)) $ is well-defined and measurable.
\end{lemma}
\begin{proof}
We  regard $\Omega\times M$ as the base space,  and define a random probability measure $(\omega,x)\mapsto\widetilde{\mu}_{(\omega,x)}:=\mu_{\omega}$  on $M$. By the measurability of $\operatorname{graph}(A)$, the function
\[
h_A:M\times(\Omega\times M)\to\mathbb R,\quad (y, (\omega,x))\mapsto \textbf{1}_{\operatorname{graph}(A)}(\omega,x,y)
\]
is bounded and measurable. Therefore, by
\cite[Proposition~3.3(i)]{MR1993844}, the map
\[
(\omega,x)\longmapsto
\int_M h_A\bigl(y,(\omega,x)\bigr)\,
\mathrm d\widetilde{\mu}_{(\omega,x)}(y)
\]
is measurable. Finally,
\[
\begin{aligned}
\int_M h_A\bigl(y,(\omega,x)\bigr)\,
\mathrm d\widetilde{\mu}_{(\omega,x)}(y)
&=
\int_M
\mathbf 1_{\operatorname{graph}(A)}(\omega,x,y)\,
\mathrm d\mu_\omega(y)\\
&=
\mu_\omega(A(\omega,x)).
\end{aligned}
\]
Thus $(\omega,x)\mapsto\mu_\omega(A(\omega,x))$ is measurable.
\end{proof}
\begin{lemma}\label{24-8-27-1044}
    For any  asymptotic  sequence $(\mathfrak{r}, A)$  on $\Omega\times M$, the Brin-Katok local  entropy  function $(\omega,x)\mapsto \overline{h}_{\mu}^{(\mathfrak{r},A)}(\varphi)(\omega,x)$ is Borel measurable.
 \end{lemma}
 \begin{proof}
    Fix a Borel  measurable function  $\mathbf{r}: \Omega\times M\to (0,1)$ and $n\in\mathbb{N}$.  First, we claim that  the graph of set-valued map 
    $$f_{\mathbf{r}}: \Omega\times M\to 2^{M},\quad (\omega,x)\mapsto B(x, \mathbf{r}(\omega,x))$$
    is a Borel measurable subset of $\Omega\times M\times M$.  Let $\mathbf{s}_{k}$ a sequence of monotone increasing simple functions on $\Omega\times M$ such that $\lim_{k\to\infty}\mathbf{s}_{k}=\mathbf{r}$.  It is obvious that 
    $$\text{graph}(f_{\mathbf{r}})=\bigcup_{k\in\mathbb{N}} \text{graph}(f_{\mathbf{s}_k}),$$
   where $f_{\mathbf{s}_k}$ is a set-valued mapping defined analogously  as $f_{\mathbf{r}}$.   The claim holds because   $\text{graph}(f_{\mathbf{s}_k})$ is  a Borel measurable subset of  $\Omega\times M\times M$ for each $k\in\mathbb{N}$.

    Denote $\mathcal{I}_n$ as the collections of all non-empty subsets of $\{0,1,\cdots,n-1\}$.   According to above claim, note that  for any  $I\in\mathcal{I}_n$ 
    $$\{(\omega,x,y)\in\Omega\times M\times M: \varphi^i_{\omega}y\in B\big(\varphi^i_{\omega}x, \mathbf{r}(\varPhi^i(\omega,x))\big)\cap A_{\theta^i\omega}\text{ for  }i\in I\}$$  is a Borel measurable subset in $\Omega\times M\times M$.  Therefore, 
    \begin{align*}
        &\big\{(\omega,x,y)\in\Omega\times M\times M: \varphi^i_{\omega}y\in B\big(\varphi^i_{\omega}x, \mathbf{r}(\varPhi^i(\omega,x))\big)\cap A_{\theta^i\omega}
        \\&\qquad\text{ for  }\varPhi^{i}(\omega,x)\in A\text{ and }i\in\{0,\dots n-1\}\big\}
        \\=&\bigcup_{I\in\mathcal{I}_n}\{(\omega,x,y)\in\Omega\times M\times M: \varphi^i_{\omega}y\in B\big(\varphi^i_{\omega}x, \mathbf{r}(\varPhi^i(\omega,x))\big)\cap A_{\theta^i\omega}\text{ for  }i\in I\}
        \\&\quad\cap\{(\omega,x,y)\in\Omega\times M\times M:  \varPhi^{i}(\omega,x)\in A\text{ for }i\in I \}
         \\&\quad\cap\{(\omega,x,y)\in\Omega\times M\times M:  \varPhi^{i}(\omega,x)\notin A\text{ for }i\in\{0,1,...,n-1\} \setminus I\}
    \end{align*} 
    is also a  Borel measurable subset of  $\Omega\times M\times M$. By \Cref{24-10-27-1252},  the mapping 
    $$(\omega,x)\mapsto\mu_{\omega}(B_{A}(n, \mathbf{r},\omega,x))$$ 
    is Borel measurable.  Hence that \Cref{24-8-27-1044} holds. 
\end{proof}

\begin{remark}
The various kinds of lower Brin--Katok local entropy are defined analogously to \Cref{24-9-6-1111} by replacing $\limsup$ with $\liminf$. The corresponding concept is denoted by replacing $\overline{\phantom{x}}$ with $\underline{\phantom{x}}$. The properties of upper Brin--Katok local entropy in this subsection also hold for lower Brin--Katok local entropy.
\end{remark}





\section{Lower bound of the measure-theoretic entropy}
\label{26-3-15-1525}
The aim of this section is to prove that the measure-theoretic entropy bounds from below the integral of the upper Brin--Katok local entropy in the non-compact random setting.
\begin{thm}
\label{24-8-28-0948}
Assume that the  $C^1$  measure-preserving  RDS $(\varphi,\mu)$ on $M$ over an invertible  metric system $(\Omega,\mathscr{F},\mathbb{P},\theta)$ 
satisfies    \eqref{2026-2-15-1654}.  Then   one has that
$$\int_{\Omega\times M}\lim_{\epsilon\to0}\limsup_{n\to\infty}-\frac{1}{n}\log\mu_{\omega}(B(\omega,x,n, \epsilon))\mathrm{d}\mu(\omega,x)\leq h_{\mu}(\varphi).$$
\end{thm}

For measure-preserving dynamical systems, the Shannon--McMillan--Breiman theorem relates metric entropy to the exponential decay of partition atoms along typical orbits. This suggests that local entropy can be compared with measure-theoretic entropy once one has a suitable measurable partition. In the non-compact setting, however, such partitions cannot be chosen with uniformly small diameter on the whole phase space. We therefore introduce a modified Brin--Katok local entropy associated with an asymptotic sequence of radii.



\subsection{Construction of  measurable partitions}

This subsection is devoted to constructing a sequence of suitable  partitions  that allow one to relate the measure-theoretic entropy of $(\varphi,\mu)$ with the Brin-Katok local entropy.



Given a  measurable partition $\alpha$ of $\Omega\times M$, recall that  we denote 
$\alpha_{\omega}:=\{A_{\omega}: A\in\alpha\}$  as the fiber  measurable partition on $M$ induced by $\alpha$ at $\omega\in\Omega$. Set 
$$\diam(\alpha_{\omega}):=\max\{\diam(B): B\in\alpha_{\omega}\}$$
when the fiber partition $\alpha_{\omega}$ is finite.  

\begin{lemma}
\label{24-8-24-1118}
Fix a  compact subset $K$  of $\Omega\times M$.  For any decreasing sequence $\{r_n\}\subset (0,1)$  converging to $0$, there  exists a sequence of non-negative decreasing random variables $\{\delta_{n}:\Omega\to \mathbb{R}_+\}$ converging   to $0$\footnote{It means that $\lim_{n\to\infty}\delta_n(\omega)=0$ for $\mathbb{P}$-a.s. $\omega\in\Omega$.},  and a sequence of countable measurable partitions $\{\alpha_n\}$ of  $K$, such that for  any $n\in\mathbb{N}$, we have 
\[\diam((\alpha_n)_{\omega})\leq 2r_n, \, \text{ and }\,  \#(\alpha_n)_{\omega}\leq r_n^{-D-\delta_{n}(\omega)},  \] 
for $\mathbb{P}$-a.s. $\omega\in \Omega$, where $D$ is the dimension of the Riemannian manifold $M$. 
\end{lemma}
\begin{proof}
By the measurable selection (for example, see \cite[Proposition 1.6.2 and Proposition 1.6.3]{A}), there exists a sequence of measurable mappings $\{\mathbf{x}_n: \Omega\to M\}_{n\in\mathbb{N}}$ such that 
\begin{align}
\label{24-8-24-943}
\overline{\{\mathbf{x}_n(\omega): n\in\mathbb{N}\}}=K(\omega),
\end{align}
for any  $\omega\in\Omega$. For notational convenience, we denote $X_J(\omega)=\{\mathbf{x}_n(\omega)\}_{n\in J}$ for any $J\subset \mathbb N$, and for any $A,B\subset M$ denote 
\[\overline{\text{dist}}(A,B): = \sup_{a\in A}\inf_{b\in B}d(a,b)\]
the (directed) Hausdorff distance from $A$ to $B$. In addition, given any  $k\in \mathbb N$ and $\epsilon>0$,  let  
\[\mathcal{I}_k=\left\{ I \subset \mathbb{N}:  \#I=k \right\},  \quad g_{k}(\omega):=\inf_{I\in\mathcal{I}_k}\overline{\text{dist}}(X_{\mathbb N}(\omega),X_{I}(\omega)),\]
and set 
\begin{align}\label{e.031301}
    \hat{\mathbf{N}}(K(\omega),\epsilon):=\inf\{k\in\mathbb{N}: g_{k}(\omega)<\epsilon \}.
\end{align}
The meaning of $g_{k}(\omega)$ is that there exists $J\in\mathcal{I}_k$ (may depend on $\omega$) such that for any $r>g_{k}(\omega)$, balls of radius $r$ centered at points from $X_J(\omega)$ cover $K(\omega)$. Then $\hat{\mathbf{N}}(K(\omega),\epsilon)$ is the minimal number of such points if the radius is fixed at $\epsilon$. Note that  $\omega\mapsto \hat{\mathbf{N}}(K(\omega),\epsilon)$ is measurable as each $g_k$ is measurable. 

\vskip0.06in 

Recall the  fractal dimension  of a subset $S$ of $M$ is defined as 
    \[\dim_{\mathrm{f}}(S):=\limsup_{\epsilon\to0}\frac{\log \mathbf{N}(S,\epsilon)}{\log(1/\epsilon)},\]
and  $\mathbf{N}(S,\epsilon)$ is the minimal  number of balls with radius $\epsilon$ centered in $S$ that are required to cover $S$. It is well-known that   $\dim_{\mathrm{f}}(S)\leq D.$ We claim that 
\begin{align}
    \label{24-8-27-1634}
    \dim_{\mathrm{f}}(K(\omega))=\limsup_{\epsilon\to0}\frac{\log \hat{\mathbf{N}}(K(\omega),\epsilon)}{\log\left(1/\epsilon\right)}.
\end{align}
Indeed, recall from  that $\mathbf{N}(K(\omega),\epsilon)$ is the minimal  number of balls with radius $\epsilon$ centered in $K(\omega)$ that are required to cover $K(\omega)$. Set $n=\mathbf{N}(K(\omega),\epsilon)$ and denote $\{y_1,\dots, y_n\}\subset K(\omega)$ as the centers such that 
$$K(\omega)\subset\bigcup_{i=1}^nB(y_i,\epsilon).$$
By density \eqref{24-8-24-943}, there exists $I\in\mathcal{I}_n$ such that 
$$K(\omega)\subset \bigcup_{i\in I}B(\mathbf{x}_i(\omega), 2\epsilon).$$ 
In view of \eqref{e.031301}, this implies that $\hat{\mathbf{N}}(K(\omega), 2\epsilon)\leq  \mathbf{N}(K(\omega),\epsilon)$. On the other hand, we have $\mathbf{N}(K(\omega),\epsilon)\leq \hat{\mathbf{N}}(K(\omega), \epsilon/2)$. Therefore, the claim \eqref{24-8-27-1634} holds. 

\vskip0.06in 

Next, for the given decreasing sequence  $\{r_m\}\subset (0,1)$ and any  $n\in\mathbb{N}$,  we let 
\[\delta_{n}(\omega)=\sup_{m\geq n}\left\{\frac{\log\hat{\mathbf{N}}(K(\omega), r_m)}{\log\left(1/r_m\right)}-\dim_{\mathrm{f}}(K(\omega)),0\right\}.\]
By \eqref{24-8-27-1634}, we know that $\lim_{n\to\infty}\delta_{n}(\omega)=0$ for $\mathbb{P}$-a.s. $\omega\in\Omega$. It also follows from the fomula of $\delta_{n}(\omega)$ that 
\begin{align}\label{e.031402}
\hat{\mathbf{N}}(K(\omega), r_n)\leq r_n^{-\dim_{\mathrm{f}}(K(\omega))-\delta_{n}(\omega)}.
\end{align}
To construct the desired partition, for any $i\in\mathbb{N}$ and $I\in\mathcal{I}_i$, we introduce the following sets, 
\begin{align}\label{e.031401}
    \begin{split}
        \Omega_i&:=\{\omega\in\Omega: \hat{\mathbf{N}}(K(\omega), r_n)=i\}, \\
        \Omega_i(I)&:=\left\{\omega\in\Omega_i:  \overline{\text{dist}}(X_{\mathbb N}(\omega), X_I(\omega))<r_n\right\}\text{ and }\\
        K_{i}(I)&:=K\cap (\Omega_i(I)\times M).
    \end{split}
\end{align}
In particular, the set $\Omega_i(I)$ consists of those samples $\omega$ sharing the same deterministic index set $I$ to cover $K(\omega)$, which is non-empty for at least some $I\in \mathcal{I}_i$.  Besides, for any $r>0$ and measurable mapping  $\mathbf{x}: \Omega\to M$, we denote 
\[\mathcal{B}(\mathbf{x}, r):=\left\{(\omega,x)\in \Omega\times M:  \omega\in\Omega\text{ and } x\in B(\mathbf{x}(\omega), r)\right\}\]
as the $r$ neighborhood along $M$ of the graph of $\mathbf{x}$. 

\vskip0.06in

Now for $I=\{j_1,j_2,\cdots,j_i\}$ and  $\omega\in \Omega_i(I)$, the balls of radius $r_n$ centered at points in $X_I(\omega)$  cover $K(\omega)$. Therefore, we can construct a finite measurable partition $\alpha^I_i=\{B_{j_1}^I,\dots, B_{j_i}^I\}$ of the block $K_{i}(I)$ as following,  
\[ B_{j_1}^I=\mathcal{B}(\mathbf{x}_{j_1}, r_n)\cap K_{i}(I),\quad B_{j_k}^I=\big(\mathcal{B}(\mathbf{x}_{j_k}, r_n)\cap K_{i}(I)\big)\setminus\bigcup_{l=1}^{k-1}B_{j_l}^I,\, k=2,\dots, i.\]
It is clear that  for any $\omega\in\Omega_i(I)$ one has that 
\begin{align}\label{e.031403}
    \begin{split}
        &\diam(\alpha_i^I(\omega))\leq 2r_n, \text{ and }\\
        &\#(\alpha_i^I(\omega))=\hat{\mathbf{N}}(K(\omega), r_n)\leq r_n^{-\dim_{\mathrm{f}}(K(\omega))-\delta_{n}(\omega)}\leq r_n^{-D-\delta_{n}(\omega)}
    \end{split}
\end{align}
by the estimate \eqref{e.031402}. As $\mathcal{I}_i$ is countable, one can write $\mathcal{I}_i=\{I_k: k\in \mathbb N\}$, where each corresponding $K_i(I_k)$ has a finite measurable partition $\alpha_i^{I_k}=\{B_{j_1}^{I_k},\dots, B_{j_i}^{I_k}\}$ satisfying \eqref{e.031403}. Defining recursively
\begin{align*}
    \widetilde{\alpha}^{I_1}=\left\{B_{j_1}^{I_1},\dots, B_{j_i}^{I_1}\right\}, \quad \widetilde{\alpha}^{I_k}=\left\{B_{j_1}^{I_k}\setminus\bigcup_{\ell=1}^{k-1}K_i(I_\ell),\dots, B_{j_i}^{I_k}\setminus\bigcup_{\ell=1}^{k-1}K_i(I_\ell)\right\}
\end{align*}
for $k\geq 2$, we obtain a countable partition $\alpha^{K_{i}}$ of the block
\[K_{i}:=K\cap (\Omega_i\times M)\]
formed by sets in $\widetilde{\alpha}^{I_k}, \, k\geq 1$. Since 
\begin{align*}
    K =\bigsqcup_{i\in\mathbb{N}} K_{i}
\end{align*}
we then obtain the desired countable measurable partition $\alpha_n$ formed by sets from the partitions $\alpha^{K_{i}}, \, i\geq 1$. 
\end{proof}

In what follows, we will use the above lemma to   construct a particular partition associated to a positive-$\mu$ measurable set $A$ and a measurable function $\mathbf{r}:\Omega\times M\to (0,1)$ that is  $\log$ $\mu$-integrable, such that the partition has a finite conditional static  entropy. This construction allows us to estimate the metric entropy from below by the Brin-Katok local entropy in the next subsection. 

\begin{lemma}
\label{24-8-27-2241}
Given   a $\mu$-positive compact  subset $K$ of $\Omega\times M$ and a small constant $\epsilon>0$, there exists a measurable subset $A\subset K$ such that  $\mu(A)\geq e^{-\epsilon}\mu(K)$ such that for any measurable function $\mathbf{r}: \Omega\times M\to (0,1)$ with $\log\mathbf{r}$  being  $\mu$-integrable, one can construct   a countable measurable partition $\alpha$ with finite conditional static entropy such that for $\mu$-a.s. $(\omega,x)\in A$, one has that 
\[\alpha(\omega, x)\subset A,\quad \text{and}\quad \diam(\alpha_{\omega}(x))\leq\mathbf{r}(\omega,x).\]
\end{lemma}
\begin{proof}
By \Cref{24-8-24-1118} with $r_n=\frac{\exp(-(n+1))}{2}$ for $n\in\mathbb{N}$,   there  exists a sequence of non-negative decreasing random variables $\{\delta_{n}:\Omega\to (0,1)\}$,  and a sequence of countable measurable partitions $\{\alpha_n\}$ of  $K$, such that for $\mathbb{P}$-a.s. $\omega\in \Omega$
\[\diam((\alpha_n)_{\omega})\leq 2r_n, \, \text{ and }\,  \#(\alpha_n)_{\omega}\leq r_n^{-D-\delta_{n}(\omega)}. \] 
 According to Egorov's theorem,   there exists $n_0\in\mathbb{N}$  and a positive $\mathbb{P}$-measurable set $\hat{\Omega}$ such that  $\delta_{n}(\omega)\leq 1$  for any $\omega\in\hat{\Omega}$ and $n\geq n_0$,  and $\mu(A)\geq e^{-\epsilon}\mu(K)$, where 
$$A:=K\cap (\hat{\Omega}\times M).$$

\vskip0.06in
 Fix a   measurable function $\mathbf{r}: \Omega\times M\to (0,1)$ with $\log\mathbf{r}$  being  $\mu$-integrable.  For $n\in\mathbb{Z}_+$, define 
\[U_n=\left\{(\omega,x)\in A: \exp(-(n+1))<\mathbf{r}(\omega,x)\leq\exp(-n)\right\}.\] 
Since $\int_{\Omega\times M}-\log\mathbf{r}(\omega,x)\mathrm{d}\mu(\omega,x)<\infty,$
it follows that 
\begin{align}
\label{24-8-27-2221}
\int_{\pi_{\Omega}(A)}\sum_{n\in\mathbb{Z}_+}n\mu_{\omega}(U_n(\omega))\mathrm{d}\mathbb{P}(\omega)<\infty.
\end{align}
By \cite[Lemma 1]{MR627789} or \cite[Lemma 4.1]{MR3038055}, one has 
\begin{align*}
&\int_{\pi_{\Omega}(A)}-\sum_{n\in\mathbb{Z}_+}\mu_{\omega}(U_n(\omega))\log\mu_{\omega}(U_n(\omega))\mathrm{d}\mathbb{P}(\omega)
\\\leq&\int_{\pi_{\Omega}(A)}\sum_{n\in\mathbb{Z}_+}n\mu_{\omega}(U_n(\omega))\mathrm{d}\mathbb{P}(\omega)+2e^{-1}\sum_{n\in\mathbb{Z}_+}e^{-n/2}<\infty.
\end{align*}

\vskip0.06in

We now show that the partition $\alpha=\{A_k: k\in\mathbb{Z}_+\}$ given below is the desired measurable countable partition, 
\begin{align*}
    A_0&=(\Omega\times M)\setminus A, \\
    \{A_k: k\in\mathbb N\} &= \left\{P\cap U_n: P\in\alpha_{{n_0+n}}\text{ and } \mu(P\cap U_n)>0,\, n\in\mathbb{Z}_+\right\}
\end{align*}
Indeed, noting that for any $n\in \mathbb N$,  $P\in \alpha_{{n_0+n}}$, and $(\omega,x)\in P\cap U_n\subset A$, one has that 
$$\diam(\alpha_{\omega}(x))\leq 2r_{n_0+n}\leq \mathbf{r}(\omega,x).$$
As $A\subset K$, this implies that $\diam(\alpha_{\omega}(x))\leq\mathbf{r}( \omega,x)$ for $\mu$-a.s. $(\omega,x)\in A$.  It remains to show the finiteness of the conditional static entropy of the partition $\alpha$. Note that \[-x\log x \leq e^{-1}, \quad x\in (0,1).\]
By \eqref{24-8-27-2221},  one has that 
\begin{align*}
    H_{\mu}(\alpha|\pi_{\Omega})&=-\int_{\Omega}\mu_{\omega}(A_0(\omega))\log \mu_{\omega}(A_0(\omega))\mathrm{d}\mathbb{P}(\omega)-\int_{\Omega}\sum_{k\in\mathbb N}\mu_{\omega}(A_k(\omega))\log \mu_{\omega}(A_k(\omega))\mathrm{d}\mathbb{P}(\omega)
    \\&\leq e^{-1}-\int_{\Omega}\sum_{n\in\mathbb{N}}\Big(\sum_{\substack{Q(\omega)\in\alpha(\omega)\\ Q(\omega)\subset U_n(\omega)}}\mu_{\omega}(Q(\omega))\log\mu_{\omega}(Q(\omega))\Big)\mathrm{d}\mathbb{P}(\omega)
    \\&\leq e^{-1}+\int_{\pi_{\Omega}(A)}\sum_{n\in\mathbb{N}}\mu_{\omega}(U_n(\omega))\Big(\log\#\alpha_{n+n_0}(\omega)-\log\mu_{\omega}(U_n(\omega))\Big)\mathrm{d}\mathbb{P}(\omega)
    \\&\leq e^{-1}+\int_{\pi_{\Omega}(A)}\sum_{n\in\mathbb{N}}\big(1+\dim_{\mathrm{f}}(K(\omega))\big)(n+1+n_0)\mu_{\omega}(U_n(\omega))\mathrm{d}\mathbb{P}(\omega)
    \\&\qquad-\int_{\pi_{\Omega}(A)}\sum_{n\in\mathbb{N}}\mu_{\omega}(U_n(\omega))\log\mu_{\omega}(U_n(\omega))\mathrm{d}\mathbb{P}(\omega)
    \\&<\infty,
\end{align*}
where in the second inequality we used the following well-known inequality  
\[-\sum_{i=1}^m x_i\log x_i\leq\sum_{i=1}^mx_i\Big(\log m-\log\sum_{i=1}^mx_i\Big) \]
valid for $x_1,\dots, x_m\in (0,1]$.  This completes the proof of \Cref{24-8-27-2241}.
\end{proof}
\subsection{Proof of \Cref{24-8-28-0948}}
The paper now proceeds as follows. We first construct measurable partitions on compact subsets of $\Omega\times M$ with controlled fiber diameters and controlled fiber cardinalities. We then refine this construction to obtain partitions of finite conditional entropy adapted to a given measurable radius function. These partitions will be used in the proof of \Cref{24-8-28-0948}. Denoting that 
\begin{align}\label{e.042201}
    \begin{split}
        &\mathbf{L}(\omega,x)=\log^+\sup_{y\in B(x,1)}|\mathrm{d}_{y}\varphi_{\omega}|,\quad A_{\mathbf{L}}(\omega,x):=\sup_{n\in\mathbb{N}}\frac{1}{n}\sum_{i=0}^{n-1}\mathbf{L}(\varPhi^i(\omega,x)),
    \end{split}
\end{align}
then $A_{\mathbf{L}}(\omega,x)$ is finite for $\mu$-a.s. $(\omega,x)\in\Omega\times M$, by Birkhoff ergodic theorem.  
\begin{proposition}
\label{25-2-9-2220}
For any  small constant $\epsilon>0$, there exists an asymptotic sequence $(\mathfrak{r}, A)$ such that   
\begin{itemize}
\item $\mu(A)\geq 1-\epsilon$;
\item for $\mu$-a.s. $(\omega,x)\in A$, one has that $\overline{h}_{\mu}^{(\mathfrak{r},A)}(\varphi)(\omega,x)\geq \overline{h}_{\mu}(\varphi)(\omega,x)$; 
\item for any $m\in\mathbb{N}$ , there is   a countable measurable partition $\alpha_m$ of  $\Omega\times M$ with finite conditional static entropy such that 
 for $\mu$-a.s. $(\omega,x)\in A$ one has that 
\begin{align}
\label{26-2-22-2108}
\alpha_m(\omega, x)\subset A,\quad \text{and}\quad \diam((\alpha_m)_{\omega}(x))\leq\mathbf{r}_m(\omega,x).
\end{align}
\end{itemize}
\end{proposition}
\begin{proof}

For  $l>0$, denote 
\begin{align}
\label{25-1-13-1518}
K_l: =\{(\omega,x)\in\Omega\times M:  A_\mathbf{L}(\omega,x)\leq l/3\}.
\end{align}
Given  $\epsilon>0$, we can choose a large  constant $l>0$  and a compact subset $K$ of $K_l$ such that  $\mu(K)\geq 1-\epsilon/2$.
Applying \Cref{24-8-27-2241} to the set $K$,  there is a measurable subset $A\subset K$ such that  $\mu(A)\geq 1-\epsilon$ such that for any measurable function $\mathbf{r}: \Omega\times M\to (0,1)$ with $\log\mathbf{r}$  being  $\mu$-integrable, one can construct   a countable measurable partition $\alpha$ with finite conditional static entropy such that for $\mu$-a.s. $(\omega,x)\in A$, one has that 
\[\alpha(\omega, x)\subset A,\quad \text{and}\quad \diam(\alpha_{\omega}(x))\leq\mathbf{r}(\omega,x).\]
Next,  we will define the asymptotic sequence  $\mathfrak{r}=\{\mathbf{r}_{m}:\Omega\times M\to\mathbb{R}_+\}_{m\in\mathbb{N}}$ with  $A$ as following 
    \begin{align}
        \label{25-11-10-1225}
        \mathbf{r}_{m}(\omega,x)=
        \begin{cases}
(me^{l})^{-\mathbf{J}(\omega,x)}&\quad \text{if } (\omega,x)\in A,
    \\ 1/m&\quad\text{if } (\omega,x)\notin A,
        \end{cases} 
    \end{align}
    where  $\mathbf{J}(\omega,x)=\inf\{i\in\mathbb{N}: \varPhi^i(\omega,x)\in A\}$ is the first return time. It remains to show that $\log\mathbf{r}_{m}$ is integrable for each $m\in\mathbb{N}$. Note that 
    \begin{align*}
    &\int_{\Omega\times M}\log\mathbf{r}_{m}(\omega,x)\mathrm{d}\mu(\omega,x)
    \\\geq&-\int_{\Omega\times M\setminus A}\log m\,\mathrm{d}\mu(\omega,x)+\int_{A}\mathbf{J}(\omega,x)(-\log m-l)\mathrm{d}\mu(\omega,x)
    \\\geq &(\mu(A)-1)\log m-(\log m+l)>-\infty,
    \end{align*}
    where we used  Kac's formula for the last  second inequality (see \cite[Theorem 2.44]{EW}). Therefore, $(\mathfrak{r},A)$ is an asymptotic sequence.
    
\begin{claim}
\label{25-2-9-2218}
Fixing   $\epsilon\in (0,1)$, $n\in\mathbb{N}$ and $(\omega,x)\in A$,   for any $y\in B(x,\epsilon e^{-nl})$ one has that 
\begin{align}
d(\varphi^i_{\omega}x,\varphi^i_{\omega}y)\leq\epsilon
\end{align}
for any $i=0,1\cdots,n$.
 \end{claim}

\begin{proof}
For any $y\in B\big(x, \frac{\epsilon}{e^{nl}}\big)$,    we will use the induction method to prove that 
$d(\varphi^i_{\omega}y,\varphi^i_{\omega}x)\leq\epsilon$ for $i=0,1,\dots,n$. Firstly, one has that 
$$d(x,y)\leq \frac{\epsilon}{e^{nl}}\leq \epsilon\ll 1 .$$
Suppose that given  $0\leq k<n$, one has that  $d(\varphi^i_{\omega}x,\varphi^i_{\omega}y)\leq \epsilon$ for any $0\leq i\leq k$.
By  \eqref{e.042201} and \eqref{25-1-13-1518} one has that
\begin{align}
\label{25-1-23-1430}
d(\varphi_{\omega}^{k+1}x, \varphi_{\omega}^{k+1}y)\leq& e^{\sum_{i=0}^{k}\mathbf{L}(\varPhi^i(\omega,x))}d(x,y)
\\\notag\leq &\epsilon e^{-nl}e^{(k+1)A_\mathbf{L}(\omega,x)}
\\\notag\leq &\epsilon e^{(k+1/3-n)l}\leq\epsilon.
\end{align}
Hence that   above claim  holds.
  \end{proof}
 For any $m\in\mathbb{N}$ and $y\in B_A(\omega,x ,n , \mathbf{r}_m)$ , one immediately has that $y\in B(\omega,x, n, 1/m)$ by \Cref{25-2-9-2218} and \eqref{25-11-10-1225}.  It follows that  $\overline{h}_{\mu}^{(\mathfrak{r},A)}(\varphi)(\omega,x)\geq \overline{h}_{\mu}(\varphi)(\omega,x)$. 
 
 Finally, applying \Cref{24-8-27-2241} for the asymptotic sequence $\mathfrak{r}=\{\mathbf{r}_m\}_{m\in\mathbb{N}}$,  we know that  for any $m\in\mathbb{N}$ ,  there is   a countable measurable partition $\alpha_m$ of  $\Omega\times M$ such that its conditional static entropy is finite 
 and  \eqref{26-2-22-2108} holds for $\mu$-a.s. $(\omega,x)\in A$. Therefore,  this completes  the all  proofs of \Cref{25-2-9-2220}.
    \end{proof}
Now, recall a version of  Shannon-McMillan-Breiman theorem in the RDS setting. 

\begin{thm}[Shannon-McMillan-Breiman  theorem for RDS]
\label{24-11-1-1933}
Let $(\varphi,\mu)$ be a   measure-preserving RDS on a complete Riemannian manifold  $M$ over an  invertible  metric system $(\Omega,\mathscr{F},\mathbb{P},\theta)$.  For any countable measurable  partition $\alpha$ of $\Omega\times M$ with $H_{\mu}(\alpha|\pi_{\Omega})<\infty$, one has that
$$h_{\mu}(\varphi,\alpha)=\int_{\Omega\times M}\lim_{n\to\infty}-\frac{1}{n}\log\mu_{\omega}\left(\big(\bigvee_{i=0}^{n-1}(\varphi^i_{\omega})^{-1}\alpha_{\theta^i\omega}\big)(x)\right)\mathrm{d}\mu(\omega,x)$$
 where  $\alpha_{\omega}$ is the fiber measurable partition on $M$ induced by $\alpha$. 
\end{thm}
\begin{proof}
Note that measure-preserving dynamical system $(\Omega\times M, \mathscr{F}\otimes\mathscr{B}( M),\mu, \varPhi)$ is not invertible.   Hence that we  can't directly obtain the above Shannon-Mcmillan-Breiman theorem  for RDS from the classical conditional Shannon-Mcmillan-Breiman theorem.    Assume that 
$$\pi: (X,\mathscr{B}(X),\mu, T)\to (Y,\mathscr{B}(Y),\nu, S)$$ be a factor map between two metric systems, and $(Y,\mathscr{B}(Y),\nu, S)$ is an invertible measure-preserving dynamical systems.  In fact, we only need to prove that  for any countable measurable  partition $\alpha$ of $X$ with $H_{\mu}(\alpha|\pi)<\infty$, one has that
$$h_{\mu}(T,\alpha|\pi)=\int_{X}\lim_{n\to\infty}-\frac{1}{n}\log\mu_{\pi(x)}\left(\big(\bigvee_{i=0}^{n-1} T^{-i}\alpha\big)(x)\right)\mathrm{d}\mu(x),$$
where $\mu=\int_Y\mu_y\,\mathrm{d}\nu(y)$ is the measure disintegration of $\mu$ with respect to $\pi$.

To use the classical conditional Shannon-Mcmillan-Breiman theorem (for example, see \cite[Theorem 3.2]{ELW} or \cite[Theorem B.0.1]{MR2809170}), we consider the invertible extension of metric  system $ (X,\mathscr{B}(X),\mu, T)$ as a factor map
$$\pi_1: (\bar{X},\mathscr{B}(\bar{X}),\bar{\mu}, \bar{T})\to (X,\mathscr{B}(X),\mu, T).$$
Then one has that 
$$\bar{\pi}:=\pi\circ \pi_1: (\bar{X},\mathscr{B}(\bar{X}),\bar{\mu}, \bar{T})\to  (Y,\mathscr{B}(Y),\nu, S)$$
is also a factor map between two invertible measure-preserving dynamical systems and $(\pi_1)_*\bar{\mu}_y=\mu_y$ for $\nu$-a.s. $y\in Y$, where $\bar{\mu}=\int_Y\bar{\mu}_y\mathrm{d}\nu(y)$ is the measure disintegration with respect to $\bar{\pi}$. Therefore,   for any countable measurable  partition $\alpha$ of $X$ with $H_{\mu}(\alpha|\pi)<\infty$, by using the conditional Shannon-Mcmillan-Breiman theorem  we know that 
$$h_{\bar{\mu}}(\bar{T},\pi_1^{-1}\alpha|\bar{\pi})=\int_{\bar{X}}\lim_{n\to\infty}-\frac{1}{n}\log\bar{\mu}_{\bar{\pi}(\bar{x})}\left(\big(\bigvee_{i=0}^{n-1} \bar{T}^{-i}\pi^{-1}_{1}\alpha\big)(\bar{x})\right)\mathrm{d}\bar{\mu}(\bar{x}).$$
It follows that 
\begin{align*}
h_{\mu}(T,\alpha|\pi)=&h_{\bar{\mu}}(\bar{T},\pi_1^{-1}\alpha|\bar{\pi})
\\=&\int_{\bar{X}}\lim_{n\to\infty}-\frac{1}{n}\log\mu_{\bar{\pi}(\bar{x})}\left(\big(\bigvee_{i=0}^{n-1} T^{-i}\alpha\big)(\pi_1(\bar{x}))\right)\mathrm{d}\bar{\mu}(\bar{x})
\\=&\int_{X}\lim_{n\to\infty}-\frac{1}{n}\log\mu_{\pi(x)}\left(\big(\bigvee_{i=0}^{n-1} T^{-i}\alpha\big)(x)\right)\mathrm{d}\mu(x).
\end{align*}
This completes the proof of \Cref{24-11-1-1933}.
\end{proof}
After completing all preparations, we begin to give the proof of \Cref{24-8-28-0948}.
\begin{proof}[Proof of \Cref{24-8-28-0948}]
Given $\epsilon>0$, let $(\mathfrak{r}, A)$ be the asymptotic sequence and $\{\alpha_m\}_{m\in\mathbb{N}}$  be a family of countable measurable partition defined as  \Cref{25-2-9-2220}, which satisfies that  
\begin{enumerate}[(1)]
\item $H_{\mu}(\alpha_{m}|\pi_{\Omega})<\infty$;
\item  for $\mu$-a.s. $(\omega,x)\in A$, one has $\alpha_{m}(\omega, x)\subset A$ and $\diam((\alpha_{m})_{\omega}(x))\leq \mathbf{r}_{m}(\omega,x)$. 
\end{enumerate}
\vskip0.06in

Fix $(\omega,x)\in A$ and $m\in\mathbb{N}$.  Note that if $y\in\big(\bigvee_{i=0}^{n-1}(\varphi^{i}_{\omega})^{-1}(\alpha_{m})_{\theta^i\omega}\big)(x)$, then for $j\in \{0,1,\cdots, n-1\}$ with  $\varPhi^{j}(\omega,x)\in A$, one has 
\begin{align*}
    &\varPhi^{j}(\omega,y)\in \alpha_{m}(\varPhi^{j}(\omega,x))\subset A, \text{ and } d(\varphi^j_{\omega}(y),\varphi^{j}_{\omega}(x))\leq \mathbf{r}_{m}(\varPhi^{j}(\omega,x)),
\end{align*}
which implies $y\in B_{A}(\omega, x, n, \mathbf{r}_{m})$. Therefore, 
by \Cref{24-11-1-1933}, for $\mu$-a.s. $(\omega,x)\in A$ one has that 
\begin{align*}
h_{\mu}(\varphi)&\geq\int_{A}\lim_{n\to\infty}-\frac{1}{n}\log\mu_{\omega}\left(\Big(\bigvee_{i=0}^{n-1}(\varphi^{i}_{\omega})^{-1}(\alpha_{m})_{{\theta^i\omega}}\Big)(x)\right)\mathrm{d}\mu(\omega,x)
\\&\geq\int_{A}\limsup_{n\to\infty}-\frac{1}{n}\log\mu_{\omega}(B_{A}(\omega, x, n, \mathbf{r}_{m}))\mathrm{d}\mu(\omega,x)
\\&\geq\int_{A}\limsup_{n\to\infty}-\frac{1}{n}\log\mu_{\omega}(B_{}(\omega, x, n, 1/m))\mathrm{d}\mu(\omega,x).
\end{align*}
It follows that $h_{\mu}(\varphi)\geq  \int _A\overline{h}_{\mu}^{}(\varphi)(\omega,x)\mathrm{d}\mu(\omega,x)$.
This finishes  the proof of \Cref{24-8-28-0948} by \Cref{25-2-9-2220} and the arbitrarily large $\mu$-measure of $A$. 
\end{proof}


\section{Upper bound of the measure-theoretic  entropy}
\label{24-4-12-1445}
In this section, we prove that the measure-theoretic entropy is bounded above by the integral of the lower Brin--Katok local entropy; see \Cref{26-2-10-2350}. The main step is to construct finite measurable partitions with null boundaries on the non-compact phase space and then adapt Katok's argument \cite{MR573822}. We also show that the
measure-theoretic entropy is finite under the integrability condition \eqref{2026-2-15-1654}. Combining \Cref{26-2-10-2350} with \Cref{24-8-28-0948}, we can obtain \Cref{26-2-13-1805}.

Given a finite measurable partition $\alpha$ of $M$, denote its boundary as 
$$\partial\alpha=\bigcup_{A\in\alpha}\partial A.$$
We say that $\alpha$ has $\nu$-null boundary if $\nu(\partial A)=0$ for every $A\in\alpha$.
\begin{lemma}
\label{26-2-10-2156}
Let $(\varphi,\mu)$ be a  measure-preserving  RDS on $M$ over an invertible metric system $(\Omega,\mathscr{F},\mathbb{P},\theta)$.  There exists a sequence of  measurable partitions $\{\beta_m\}$ of $M$ with $\mu_M$-null boundary such that 
$$h_{\mu}(\varphi)=\lim_{m\to\infty} h_{\mu}\big(\varphi,\pi_M^{-1}\beta_m\big),$$ where $\pi_M: \Omega\times M\to M$ is the coordinate projection and $\mu_M:=(\pi_M)_*\mu$.
 \end{lemma}
 \begin{proof}
Let $\{K_n\}_{n\in\mathbb N}$ be an increasing sequence of compact
subsets of $M$ such that $M=\cup_{n\in\mathbb N}K_n.$
We first construct a sequence of finite Borel measurable partitions
$\{\alpha_n\}_{n\in\mathbb N}$ with $\mu_M$-null boundary. Set $B_0=\varnothing$. Suppose that $B_{n-1}$ has already been
constructed and is relatively compact. The set
\[
K_n\cup\overline{B_{n-1}}
\]
is compact. For every point of this set, one can choose a relatively
compact open ball of diameter less than $1/n$ whose boundary has zero
$\mu_M$-measure.  By compactness of $K_n\cup \overline{B_{n-1}}$, there exist finitely many such open balls
\[
U_{n,1},\ldots,U_{n,k(n)}
\]
covering $K_n\cup\overline{B_{n-1}}$. Define
\[
C_{n,1}:=U_{n,1},
\qquad
C_{n,j}
:=
U_{n,j}\setminus\bigcup_{\ell=1}^{j-1}U_{n,\ell},
\quad 2\leq j\leq k(n),
\]
and
\[
B_n:=\bigcup_{j=1}^{k(n)}U_{n,j}
    =\bigsqcup_{j=1}^{k(n)}C_{n,j}.
\]
It follows that  the sets $C_{n,j}$ are pairwise
disjoint, relatively compact Borel sets satisfying $\diam(C_{n,j})<1/n.$
Moreover,
$$\mu_M(\partial C_{n,j})\leq \sum_{\ell=1}^j\mu_M(\partial U_{n,\ell})=0,\quad \mu_M(\partial B_n)\leq \sum_{\ell=1}^{k(n)}\mu_M(\partial U_{n,\ell})=0 .$$
Consequently,
\[
\alpha_n
:=
\left\{
C_{n,1},\ldots,C_{n,k(n)},M\setminus B_n
\right\}
\]
is a finite Borel measurable partition of $M$ with
$\mu_M$-null boundary.

By above  construction, we know that  $K_n\cup B_{n-1}\subset B_n.$
Therefore, $\{B_n\}$ is increasing and $M=\bigcup_{n\in\mathbb N}B_n.$  Now, we set $\beta_n:=\bigvee_{i=1}^{n}\alpha_i.$
Then $\{\beta_n\}_{n\in\mathbb N}$ is an increasing sequence of finite
Borel measurable partitions. Since every atom of $\beta_n$ is a finite
intersection of atoms of $\alpha_1,\ldots,\alpha_n$, each $\beta_n$
has $\mu_M$-null boundary. Let $\mathscr G:=\bigvee_{n\geq1}\widehat{\beta}_n,$
where $\widehat{\beta}_n$ denotes the $\sigma$-algebra generated by
$\beta_n$. We claim that $\mathscr G=\mathscr B(M).$
It suffices to show that every open subset of $M$ belongs to
$\mathscr G$.

Let $U\subset M$ be open and let $x\in U$. Choose $N$  sufficiently large such that
$x\in B_N$. Since $B_N\subset B_n$ for every $n\geq N$, the atom
$\alpha_n(x)$ has diameter less than $1/n$. Since $\beta_n$ refines
$\alpha_n$, we also have $\operatorname{diam}\bigl(\beta_n(x)\bigr)<1/n$
for every $n\geq N$. Choosing $n$ sufficiently large, we obtain $\beta_n(x)\subset U.$
It follows that
\[
U
=
\bigcup_{n\geq1}
\bigcup_{\substack{D\in\beta_n\\D\subset U}}D.
\]
The right-hand side is a countable union of elements of $\mathscr G$.
Thus   $\mathscr G=\mathscr B(M).$  Particularly,  we have 
\[
\pi_M^{-1}\big(\mathscr G\big)
\vee
\pi_\Omega^{-1}\big(\mathscr F\big)
=
\mathscr F\otimes\mathscr B(M).
\]
Hence that \Cref{26-2-10-2156} follows from \Cref{26-2-22-2151}.
\end{proof}
Next, we recall a technical lemma. The reader can refer to the details in \cite{MR573822}. Let $\alpha$ be a finite Borel measurable partition of $M$ and $\alpha_{\omega}^n:=\bigvee_{i=0}^{n-1}(\varphi_{\omega}^{i})^{-1}\alpha$ for  $\omega\in\Omega$ and $n\in\mathbb{N}$. Now,   we  define a family of  pseudo-metrics on $M$ induced by $\alpha$ and the RDS $\varphi$ as follows 

$$\rho_{\alpha^n_{\omega}}(x,y)=\frac{\sum_{i=0}^{n-1}\rho_{\alpha}(\varphi^i_{\omega}x,\varphi^i_{\omega}y)}{n},$$
where 
$$\rho_{\alpha}(x,y):=
\begin{cases}
0 & \text{ if }  x,y \text{ belong to the same atom of  } \alpha,
\\1 &\text{ otherwise}.
\end{cases}
$$ 
The following combinatorial estimate follows from Stirling's formula.
\begin{lemma}
\label{26-2-10-2333}
Let $(\varphi,\mu)$ be a  measure-preserving  RDS on $M$ over a metric system $(\Omega,\mathscr{F},\mathbb{P},\theta)$ and $\alpha$ be a finite measurable partition on $M$.  Given a sufficiently   small constant $r$,  there exists a large number  $N$ (depending only  on $r$ and $\alpha$) such that for any $(\omega,x)\in\Omega\times M$ and any $n\geq N$
$$\#\{ E\in\alpha^n_{\omega} :  E\cap B_{\rho_{\alpha^n_\omega}}(x,r)\neq\varnothing\}\leq  e^{ 3n \sqrt{r}\#\alpha},$$
where   $B_{\rho_{\alpha^n_\omega}}(x,r):=\{y\in M:  \rho_{\alpha^n_\omega}(x,y)\leq r\}$.
\end{lemma}

We now state and prove the main result of this section.
\begin{proposition}
\label{26-2-10-2350}
Let $(\varphi,\mu)$ be a   measure-preserving  RDS on $M$ over an invertible metric system $(\Omega,\mathscr{F},\mathbb{P},\theta)$.  
One has that 
\begin{align}
\label{25-12-23-1623}
h_{\mu}(\varphi)\leq \int_{\Omega\times M}\lim_{\epsilon\to 0}\liminf_{n\to\infty}-\frac{\log\mu_{\omega}(B(\omega,x, n, \epsilon))}{n}\mathrm{d}\mu(\omega,x).
\end{align}
 \end{proposition}
 \begin{proof}
Assume, to the contrary, that \eqref{25-12-23-1623} does not hold.  By \Cref{26-2-10-2156}, \Cref{24-11-1-1933} and Lusin's  lemma,  we can  choose   a finite  measurable partition $\alpha$ of $M$ with $\mu_M$-null boundary,  two positive constants $h\gg\delta_0$, and a  positive-$\mu$ measurable  subset  $\Lambda_1$ such that   for any $(\omega,x)\in \Lambda_1$  one has that 
\begin{align}
\lim_{n\to\infty}-\frac{\log\mu_{\omega}(\alpha_{\omega}^n(x))}{n}&>h-\delta_0/5
\\\notag&>h-\delta_0>\lim_{\epsilon\to 0}\liminf_{n\to\infty}-\frac{\log\mu_{\omega}(B(\omega,x, n, \epsilon/3))}{n}.
\end{align}

Fix a small constant $\delta>0$ such that $\delta\ll \mu(\Lambda_1)$.   Since $\partial\alpha =\lim_{\epsilon\to 0}\mathscr{U}_{\epsilon}$, where  
\begin{align}
\mathscr{U}_{\epsilon}:=\{x\in M:  B(x, \epsilon)\text{ is not contained in }\alpha(x)\},
\end{align}
there exists  a small constant $\epsilon_0$   such that for any positive constant $\epsilon\leq \epsilon_0$ one has that   $\mu_M(\mathscr{U}_{\epsilon}(\alpha))\leq \delta^2/2.$
Moreover, by the last inequality in the definition of $\Lambda_1$, after replacing
$\Lambda_1$ by a measurable subset of positive $\mu$-measure if necessary, we may choose
$0<\epsilon\leq\epsilon_0$ such that for every $(\omega,x)\in\Lambda_1$,
$$\liminf_{n\to\infty}
-\frac{\log\mu_\omega(B(\omega,x,n,\epsilon/3))}{n}<h-\frac{3\delta_0}{5}.$$

Fix the above $\epsilon\in(0,\epsilon_0)$.  By using Birkhoff ergodic theorem, we know that  
\begin{align*}
\int_{\Omega\times M}\lim_{k\to\infty}\frac{1}{k}\sum_{i=0}^{k-1}\mathbf{1}_{\pi^{-1}_M\mathscr{U}_{\epsilon}}(\varPhi^i(\omega,x))\mathrm{d}\mu(\omega,x)=\mu (\pi^{-1}_M\mathscr{U}_{\epsilon})\leq \delta^2/2.
\end{align*}
By  Egorov's theorem and $\mu(\Lambda_1)\gg\delta$,  there exists a large number $N_0$  and positive-$\mu$ measurable subset $\Lambda_2\subset \Lambda_1$  such that  for any $\omega\in \pi_{\Omega}(\Lambda_2)$ and $n\geq N_0$ one has that 
       \begin{align}
      \Lambda_2(\omega)\subset \Big\{y\in \Lambda_1(\omega): \frac{1}{k}\sum_{i=0}^{k-1}\mathbf{1}_{\mathscr{U}_{\epsilon}(\alpha)}(\varphi^i_{\omega}y)\leq\delta \text{ for any } k\geq n\Big\}.
      \end{align} 
  Therefore,  one has  that  for any $(\omega,x)\in \Lambda_2$ and $n\geq N_0$
  \begin{align}
  \label{25-2-11-1825}
  \rho_{\alpha^n_\omega}(x, y)\leq \delta,
 \end{align}
  whenever $y\in B(\omega,x, n,\epsilon)$.
 By Egorov's theorem,  there exists a positive-$\mu$ measurable subset $\Lambda_3\subset\Lambda_2$ and any sufficiently large $N_1\geq N_0$ such that for any $(\omega,x)\in \Lambda_3$ one has that  for any  $n\geq N_1$
 \begin{align}
 \label{25-2-11-1836}
 \mu_{\omega}(\alpha^n_{\omega}(x))\leq e^{-n(h-2\delta_0/5)}.
 \end{align}

Given $ N\geq N_1$ and  $\omega\in\pi_{\Omega}(\Lambda_3)$,   denote  
$$\Lambda^{n,N}_3(\omega)=\Big\{x\in\Lambda_3(\omega):  n \text{ is the smallest  integer s.t. } \mu_{\omega}(B(\omega,x,n,\epsilon/3))\geq e^{-n(h-3\delta_0/5)}\Big\},$$
for $n\geq N$.
It is clear that  $\Lambda_3(\omega)=\bigcup_{n\geq N} \Lambda^{n,N}_3(\omega)$.  For any $n\geq N$,  choose a maximal $(n,\epsilon)$-separated set
$\{x_i:1\le i\le s_n\}$ in $\Lambda_3^{n,N}(\omega)$.
Then the balls $B(\omega,x_i,n,\epsilon)$ cover
$\Lambda_3^{n,N}(\omega)$, while the balls
$B(\omega,x_i,n,\epsilon/3)$ are pairwise disjoint.
For any $x_i\in\Lambda_3^{n,N}(\omega)$, we have
$$\mu_\omega(B(\omega,x_i,n,\epsilon/3))\geq e^{-n(h-3\delta_0/5)},$$
which implies $s_n\le e^{n(h-3\delta_0/5)}$.  By \Cref{26-2-10-2333} and   \eqref{25-2-11-1825}, we have that for any $x\in\Lambda_3(\omega)$ 
$$\#\{E\in \alpha_{\omega}^n: E\cap B(\omega, x, n,\epsilon)\neq\varnothing \}\leq  e^{3n\sqrt{\delta}\#\alpha}.$$
 Therefore, there exists at most  $e^{n(h-3\delta_0/5)+3n\sqrt{\delta}\#\alpha}$-atoms in  $\alpha^n_{\omega}$ to cover  $\Lambda^{n,N}_3(\omega)$.  By  \eqref{25-2-11-1836}, one has that 
\begin{align*}
\mu_{\omega}(\Lambda_3(\omega))&=\sum_{n\geq N}\mu_{\omega}(\Lambda_3^{n,N}(\omega))
\\&\leq \sum_{n\geq N}e^{n(h-3\delta_0/5)}e^{3n\sqrt{\delta}\#\alpha} e^{-n(h-2\delta_0/5)}
\\&\leq \sum_{n\geq N}e^{3n\sqrt{\delta}\#\alpha-n\delta_0/5}.
\end{align*}
This  yields a contradiction  by choosing $\delta\ll \delta_0$ and  $N$  big enough.  Therefore,  \Cref{26-2-10-2350} holds.
 \end{proof}

Finally,  we prove that the measure-theoretic entropy is  finite under the setting of \Cref{26-2-13-1805}. 

\begin{proposition}
\label{26-7-28-0958}
Let $(\varphi,\mu)$ be a $C^1$ measure-preserving RDS on the 
Riemannian manifold $M$  over an invertible metric system
$(\Omega,\mathscr{F},\mathbb{P},\theta)$. If it satisfies 
\eqref{2026-2-15-1654}, then
\begin{align}
\label{26-7-27-entropy-bound}
h_{\mu}(\varphi)
\leq\int_{\Omega\times M}
\overline{h}_{\mu}(\varphi)(\omega,x)\,\mathrm{d}\mu(\omega,x)
<+\infty.
\end{align}
\end{proposition}
\begin{proof} 
 Set
\begin{align}
\label{26-7-27-def-L}
\mathbf{L}(\omega,x)
:=
\log^+\sup_{y\in B(x,1)}
\left|\mathrm{d}_y\varphi_\omega\right|
\quad\text{ and }\quad S_n\mathbf{L}(\omega,x):=
\sum_{i=0}^{n-1}
\mathbf{L}\bigl(\varPhi^i(\omega,x)\bigr) \text{ for } n\in\mathbb{N}.
\end{align}
By the Birkhoff ergodic theorem,  one has that for $\mu$-a.s. $(\omega,x)\in\Omega\times M$
$$\mathbf{L}^*(\omega,x)=\lim_{n\to\infty}\frac{S_n\mathbf{L}(\omega,x)}{n}$$
and $\int_{\Omega\times M}\mathbf{L}^*\,\mathrm{d}\mu=\int_{\Omega\times M}\mathbf{L}\,\mathrm{d}\mu.$ Fix such a point $(\omega,x)$,  and a constant $\epsilon\in (0,1)$.  By using a similar argument as \Cref{25-2-9-2218}, one has that 
\begin{align}
\label{26-7-27-ball-inclusion}
B\big(x,\epsilon e^{-S_n\mathbf{L}(\omega,x)}\big)\subset B(\omega,x,n,\epsilon).
\end{align}
Fix $\eta>0$. For all sufficiently large $n$, one has $S_n\mathbf{L}(\omega,x)
\leq
n\bigl(\mathbf{L}^*(\omega,x)+\eta\bigr).$
By \eqref{26-7-27-ball-inclusion}, we know that
\begin{align*}
\label{26-7-27-smaller-ball}
B\big(x,\epsilon e^{-n(\mathbf{L}^*(\omega,x)+\eta)}\big)
\subset B(\omega,x,n,\epsilon).
\end{align*}
It follows that 
\begin{align*}
&\limsup_{n\to\infty}
-\frac{1}{n}
\log\mu_\omega\bigl(B(\omega,x,n,\epsilon)\bigr)
\\
&\quad\leq
\limsup_{n\to\infty}
-\frac{1}{n}
\log\mu_\omega
\left(
B\left(x,\epsilon e^{-n(\mathbf{L}^*(\omega,x)+\eta)}\right)
\right)
\\
&\quad\leq
D\bigl(\mathbf{L}^*(\omega,x)+\eta\bigr)
\end{align*}
for $\mu$-a.s. $(\omega,x)$, where $D$ is the dimension of the manifold $M$. Therefore, by \Cref{26-2-10-2350}, we obtain that 
\begin{align*}
h_\mu(\varphi)
&\leq
\int_{\Omega\times M}
\underline{h}_\mu(\varphi)(\omega,x)\,
\mathrm{d}\mu(\omega,x)
\\
&\leq
D\int_{\Omega\times M}
\mathbf{L}^*(\omega,x)+\eta\,\mathrm{d}\mu(\omega,x)
\\
&=
D\int_{\Omega\times M}
\mathbf{L}(\omega,x)\,\mathrm{d}\mu(\omega,x)+D\eta
<+\infty.
\end{align*}
 This proves \Cref{26-7-28-0958}.
\end{proof}






\section{Margulis-Ruelle inequality}
\label{26-3-15-1526}

This section is devoted to proving that the upper neutralized Brin--Katok local entropy is bounded above by the sum of the positive Lyapunov exponents. Combined with the Brin--Katok local entropy formula in \Cref{26-2-13-1805}, this yields the Margulis--Ruelle inequality in \Cref{24-4-12-1500}.
\begin{thm}
\label{24-4-5-2040}
Assume that the  $C^2$   measure-preserving RDS $(\varphi,\mu)$ on $M$ over $(\Omega,\mathscr{F},\mathbb{P},\theta)$
satisfies {\rm \hyperref[23-4-29-1706]{(H1)}-\hyperref[24-9-26-1718]{(H2)}}. We have
\begin{align}
\int_{\Omega\times M}\lim_{\epsilon\to 0}\limsup_{n\to\infty}\frac{-\log\mu_{\omega}(B(\omega,x, n, e^{-n\epsilon}))}{n}\mathrm{d}\mu(\omega,x)\leq\int_{\Omega\times M}\lambda_{\sum}^+(\omega,x)\mathrm{d}\mu(\omega,x),
\end{align}
 where $\lambda_{\sum}^+$  is  sum of all positive Lyapunov exponents counted with multiplicities.  
\end{thm}


\medskip

The structure of the section goes as follows.  In \Cref{26-2-23-1409}, we mainly review the multiplicative ergodic theorem and unconventional  Lyapunov charts from \cite{MR1369243}.    In \Cref{s.010602}, we extend the idea of unstable graphs with controlled dispersion in \cite{MR627789} to decompose the neutralized Bowen balls. In \Cref{24-10-29-2033}, we establish the key estimate for the decaying rate of the volume measure of the neutralized Bowen ball through volume expansion rate along the unstable directions. The proof of \Cref{24-4-5-2040} is provided in \Cref{26-2-23-1418}.




\subsection{Multiplicative ergodic theorem and Lyapunov charts}
\label{26-2-23-1409}
Recall that $\varPhi$ denotes the skew product transformation on $\Omega\times M$ induced by the RDS $\varphi$ (see   \eqref{23-1-3-1520}), $D$ is the dimension of the Riemannian manifold $M$.  Applying \cite[Theorem 3.2 in Chapter 1]{MR1369243} to    the derivative  linear cocycle of the smooth measure-preserving  RDS $(\varphi,\mu)$  generated  by 
$$(\omega,x)\mapsto\mathrm{d}_x\varphi_{\omega},$$
one has the following convenient  multiplicative ergodic theorem.
\begin{thm}
\label{24-3-28-2010}
There  is  an $\varPhi$-invariant Borel subset  $\Gamma$ of $\Omega\times M$  with $\mu(\Gamma)=1$,  a measurable function  $\mathbf{r}:\Gamma\to\mathbb{N}$,  and  on the level sets of $\mathbf{r}$ a collection measurable functions $\lambda_1,\dots,\lambda_{\mathbf{r}(\omega,x)}: \Gamma\to\mathbb{R}\cup\{-\infty\}$ with 
$$\lambda_1(\omega,x)>\lambda_2(\omega,x)>\cdots>\lambda_{\mathbf{r}(\omega,x)}(\omega,x)\geq-\infty$$
and  a measurable filtration 
 $$T_xM=V^{1}(\omega,x)\supset\cdots\supset V^{\mathbf{r}}(\omega,x)\supset \{0\} =:V^{\mathbf{r}+1}(\omega,x)$$
 with following properties 
 
  \begin{enumerate}[\rm (I)]
\item for each $j=1,2,\dots, r(\omega, x)$,  $m_j(\omega,x):=\dim V^j(\omega,x)- \dim V^{j+1}(\omega,x)$ is $\varPhi$-invariant; 
  \item  Invariance:   $\mathrm{d}_x\varphi_{\omega}V^{j}(\omega,x )\subset V^{j}(\varPhi(\omega, x))$;
  \item Lyapunov exponent: 
  $$\lim_{n\to\infty}\frac{\log|\mathrm{d}_x\varphi^n_{\omega}v|}{n}=\lambda_j(\omega,x)$$ 
  for all $ v\in V^j(\omega,x)\setminus  V^{j+1}(\omega,x)$;
     \item   Temperedness: let $\{v_1,\cdots, v_D\}$ be a family of basis of $T_xM$ which satisfies 
  $$\lim_{n\to\infty}\frac{1}{n}\log|\mathrm{d}_x\varphi^n_{\omega}(v_i)|=\lambda_{(i)}(\omega,x)$$
  where $\lambda_{(i)}(\omega,x)$ is the Lyapunov exponents counted with multiplicity.  Furthermore,  we have 
    \begin{align*}
\lim_{n\to\infty}\frac{1}{n}\log\angle (V^j(\theta^n\omega,\varphi_{\omega}^nx), \mathrm{d}_x\varphi^n_{\omega} (V^j(\omega,x)^{\perp}))=0,
  \end{align*} 
  where $E^{\perp}$ is the  orthogonal complement space  to  a subspace $E$ of $T_xM$. 
  \end{enumerate}
\end{thm}

\begin{remark}
The subset $\Gamma$ and $\{(\lambda_j, m_j): 1\leq j\leq r\}$ in \Cref{24-3-28-2010} are  called as the \emph{Lyapunov regular set} and  \emph{Lyapunov spectra} of  measure-preserving RDS $(\varphi,\mu)$, respectively.
\end{remark}

Next,  we mainly review the process in \cite[Section 1 in Chapter III]{MR1369243} to construct a family of norms on  the  Lyapunov regular set.   Letting  $\mathbf{l}: \Omega\times M\times  \mathbb{Z}_+\to\mathbb{R}_+$ be a positive measurable function and $\tau$ be a small positive constant, then $\mathbf{l}$ is called \emph{$\tau$-slowly varying} if for any $(\omega,x,n)\in\Omega\times M\times\mathbb{Z}_+$ and $m\in\mathbb{Z}_+$ one has that 
\begin{align}
\mathbf{l}(\omega,x, n+m)\leq \mathbf{l}(\omega,x,n)e^{\tau m}.
\end{align}

To deal with the difficulty arising from non-ergodicity,  we need to cut the Lyapunov regular set $\Gamma$ defined as \Cref{24-3-28-2010} into several disjoint $\varPhi$-invariant measurable subsets. These invariant  subsets are called  Pesin's blocks.   Now, we give a specific form:
  
  \medskip
  For any $\text{u}\in\mathbb{N},\vec{m}=(m_1,\dots, m_\text{u})\in\mathbb{N}^\text{u}$, a vector $\vec{\lambda}=(\lambda_1,\dots,\lambda_\text{u})\in\mathbb{R}^\text{u}_+$ and a positive small constant $\tau$,  
  \begin{align}
\label{25-11-3-1637}
\Gamma(\text{u},\vec{m},\vec{\lambda}; \tau):=\Big\{&(\omega,x)\in\Gamma:  \text{there are $\text{u}$ positive Lyapunov exponents}, 
\\\notag& m_j(\omega,x)=m_j, 0< \lambda_j(\omega,x)-\lambda_j\leq\tau, j=1,\cdots, \text{u}\Big\}.
\end{align}
It is obvious that  $\Gamma(\text{u},\vec{m},\vec{\lambda}; \tau)$ is a $\varPhi$-invariant measurable subset of $\Gamma$ and 
\begin{align}
\label{25-11-12-1407}
\Gamma=\bigcup_{\text{u}\in\mathbb{N},\vec{m}\in\mathbb{N}^\text{u}, \vec{\lambda}\in\mathbb{Q}^\text{u}_+}\Gamma(\text{u},\vec{m},\vec{\lambda}; \tau),
\end{align}
for any $\tau>0$.  In most arguments  of this  paper, we will limit our version on a fixed Pesin's block.  Specifically, from here to the end of  this Section,  we fix $\text{u},\vec{m}, \vec{\lambda},\tau$, and write
$$\mathbf{\Gamma}=\Gamma(r,\vec{m},\vec{\lambda}; \tau).$$
Now,  we begin to  introduce some convenient notations. We define the unstable space\footnote{In the absence of positive Lyapunov exponents, we set $E_0^u(\omega, x)=\{0\}$.}
 and center-stable space as following:
\begin{align*}
\begin{split}
&E_0^{cs}(\omega,x):=V^{\text{u}+1}(\omega,x), \quad E_0^u(\omega,x):=\big(E_0^{cs}(\omega,x)\big)^{\perp};
\\&E_n^{cs}(\omega,x):=E_0^{cs}\big(\varPhi^n(\omega,x)\big),  \quad E_n^u(\omega,x):=\mathrm{d}_x\varphi^n_{\omega}\big(E_0^{u}(\omega,x)\big) \text{ for } n\in\mathbb{N}.
\end{split}
\end{align*}
From  now on, we denote  $\pi^{\sigma}_{(\omega,x),n}$ as  the projections onto $E_n^{\sigma}(\omega,x)$ along  $E_n^{\sigma'}(\omega,x)$ for $\sigma\neq\sigma'\in\{u,cs\}$.  Fixing a  sufficiently small constant  $\tau$, we assume that 
\begin{align}\label{25-2-9-1638}
   \chi:=\min_{}\{\lambda_i: \lambda_i>0\}-\tau\geq 99\tau>0,
\end{align} then one has that the following lemma.
\begin{lemma}[{\cite[Lemma 1.1 in Chapter III]{MR1369243}}]
\label{24-1-21-1935}
There exists a $\tau/3$-slowly varying function $\mathbf{l}_1:\Gamma\times\mathbb{Z}_+\to (0,+\infty)$ such that for any $(\omega,x,n)\in\mathbf{\Gamma}\times\mathbb{Z}_+$ and $m\in\mathbb{Z}_+$ we have that 
\begin{enumerate}[(1)]
\item  $|\mathrm{d}_{x_n}\varphi^m_{\omega_n}v|\leq \mathbf{l}_1(\omega,x,n ) e^{\tau m/3}|v|$ for any $v\in E_n^{cs}(\omega,x)$;
\item $|\mathrm{d}_{x_n}\varphi^m_{\omega_n}u|\geq \mathbf{l}_1(\omega,x,n )^{-1} e^{(\chi+\tau/3)m}|u|$ for any $u\in E_n^{u}(\omega,x)$;
\item $\angle (E^u_{n+m}(\omega,x), E^{cs}_{n+m}(\omega,x))\geq \mathbf{l}_1(\omega,x,n)^{-1}e^{-\tau m/3}$,
\end{enumerate}
where $\omega_n:=\theta^n\omega$ and $x_n:=\varphi^n_{\omega}x$.
\end{lemma}
According to above lemma, we  can define a family of  $\tau$-Lyapunov inner products   as following 
\begin{align*}
\langle w_1, w_2\rangle_{(\omega,x),n}=
\begin{cases}
\sum_{i=0}^n\langle \mathrm{d}_{x_n}\varphi^{-i}_{\omega_n}(w_1), \mathrm{d}_{x_n}\varphi^{-i}_{\omega_n}(w_2)\rangle e^{2i\chi},&\quad \text{ if } w_1,w_2\in E_n^{u}(\omega,x)
\\\sum_{i\in\mathbb{Z}_+}\langle \mathrm{d}_{x_n}\varphi^i_{\omega_n}(w_1), \mathrm{d}_{x_n}\varphi^i_{\omega_n}(w_2) \rangle e^{-2i\tau}, &\quad \text{ if } w_1,w_2\in E_n^{cs}(\omega,x)
\\ 0, &\quad\text{ otherwise},
\end{cases}
\end{align*}
where 
$$\mathrm{d}_{x_n}\varphi^{-i}_{\omega_n}=(\mathrm{d}_{x_{n-i}}\varphi_{\omega_{n-i}})^{-1} |_{E_{n-i+1}^{u}(\omega,x)} \circ \cdots \circ (\mathrm{d}_{x_{n-1}} \varphi_{\omega_{n-1}})^{-1} |_{E_n^{u}(\omega,x)}.
$$
Denote the norms as $|\cdot|_{(\omega,x),n}$ induced by the inner product $\langle\cdot,\cdot\rangle_{(\omega,x),n}$.  This family of measurable norms has  following properties.
\begin{lemma}[{\cite[Lemma 1.3 in Chapter III]{MR1369243}}]
\label{23-12-18-21-37}
\begin{enumerate}[(a)] 
\item For any $m\in\mathbb{N}$ and $(\omega,x,n)\in\mathbf{\Gamma}\times\mathbb{Z}_+$, one has that 
\begin{align*}
&|\mathrm{d}_{x_n}\varphi^{m}_{\omega_n}(\pi^u_{(\omega,x),n}w)|_{(\omega,x),n+m}\geq e^{m\chi}|\pi^u_{(\omega,x),n}w|_{(\omega,x),n},
\\&|\mathrm{d}_{x_n}\varphi^m_{\omega_n}(\pi^{cs}_{(\omega,x),n}w)|_{(\omega,x),n+m}\leq e^{m\tau}|\pi^{cs}_{(\omega,x),n}w|_{(\omega,x),n},
\end{align*}
holds for any $w\in T_{x_n}M$;
\item the mapping  $(\omega,x)\mapsto |w|_{(\omega,x)}$ is measurable for any fixed $w\in M$;
\item there exists a  $2\tau/3$-tempered function  $\mathbf{l}_2: \mathbf{\Gamma}\times\mathbb{Z}_+\to[1,+\infty)$  such that  for any $w\in T_{x_n}M$ and $n\in\mathbb{Z}_+$ one has that 
$$\frac{1}{\sqrt{2}}|w|\leq|w|_{(\omega,x),n}\leq \mathbf{l}_2(\omega,x,n)|w|.$$
\end{enumerate}
\end{lemma}

Finally, we will construct the Lyapunov chart by using the  integrability condition  \eqref{23-4-29-1706}.   It is clear that there  exist two  orthogonal subspaces $E^u$ and $E^{cs}$  of a  Euclidean space  $\mathbf{R}^D$ with  the  norm $\|\cdot\|$ induced by its inner product  such that 
\begin{itemize}
\item $\dim E^u=\dim E_n^u(\omega, x)$, and $\dim E_n^{cs}=\dim E_n^{cs}(\omega, x)$ for any $(\omega,x, n)\in\mathbf{\Gamma}\times\mathbb{Z}_+$;
\item By the measurable selection theorem (see \cite{MR467310}), there exist a family of  invertible bounded  measurable linear maps 
\begin{align}\label{e.051901}
    \{L_{(\omega,x),n} :  (T_{x_n}M,|\cdot|)\}_{(\omega,x)\in\Gamma}\to (\mathbf{R}^D,\| \cdot\|)
\end{align}
such that 
\begin{enumerate}[(i)]
\item $L_{(\omega,x),n}(E_n^{\sigma}(\omega,x))=E^{\sigma}$, where $\sigma\in\{u,cs\}$;
\item $\| L_{(\omega,x),n}w\|=|w|_{(\omega,x),n}$, wherever $w\in\mathbf{R}^D$.
\end{enumerate}
\end{itemize}

In the following context, we  always assume that the regular radius $\rho_b$  satisfies the integrable condition \eqref{24-9-26-1718} for some positive constant $b$. Now,  we introduce the connection maps along the orbit
\begin{align}
\label{24-10-28-1423}
\varphi_{(\omega,x),n}:= \mathcal{T}^{-1}_{(\omega,x), n+1}\circ \varphi_{\omega_n}\circ \mathcal{T}_{(\omega,x),n}
\end{align}
where $\mathcal{T}_{(\omega,x),n}:\mathbf{R}^D\to M $ with $w\mapsto \exp_{x_n}\circ L^{-1}_{(\omega,x),n}(w)$.   

\medskip

We now establish several properties of the connecting maps.  For convenience, we denote $\pi^{\sigma}$ as the projection on $E^{\sigma}$ along $E^{\sigma'}$, where $\sigma\neq\sigma'\in\{u,cs\}$, and 
\begin{align}\label{24-4-20-1940}
 &\mathbf{B}^{\sigma}(r):=\{w\in E^{\sigma}:  \|w\|<r\},
 \\&\mathbf{B}(r):=\{w\in\mathbf{R}^D: \|\pi^{\sigma}w\|< r\text{ for } \sigma=u,cs\}.        
\end{align}
\begin{proposition}[Lyapunov charts]
\label{24-3-8-1417}
There exists a $\tau$-slowly varying   function $\mathbf{l}: \mathbf{\Gamma}\times\mathbb{Z}_+\to[1,+\infty)$  such that for any   $\delta\in (0,1)$,  $\varphi_{(\omega,x),n}$ is well-defined on $\mathbf{B}(\delta \mathbf{l}(\omega,x,n)^{-1})$ such that the followings hold for any $(\omega,x,n)\in\Gamma\times\mathbb{Z}_+$:
\begin{enumerate}[(a)]
    \item\label{23-12-18-2136-1}  For any $w\in\mathbf{R}^D$, one has that 
    $$\mathbf{l}(\omega,x,n)^{-1}\|w\|\leq |L^{-1}_{(\omega,x),n}(w)|\leq\sqrt{2}\|w\|.$$
    \item\label{23-12-18-2136-2}  We have $\mathrm{d}_0\varphi_{(\omega,x),n}E^{u}=E^{u}$, $\mathrm{d}_0\varphi_{(\omega,x),n}E^{cs}\subset E^{cs}$, and for any $w\in\mathbf{R}^D$, 
    \begin{align*}
        &\|\mathrm{d}_0\varphi_{(\omega,x),n}(\pi^uw)\|\geq e^{\chi}\|\pi^uw\|,\quad \|\mathrm{d}_0\varphi_{(\omega,x),n}(\pi^{cs}w)\|\leq e^{\tau}\|\pi^{cs}w\|.
    \end{align*}
    \item\label{23-12-18-2136-3} Over $\mathbf{B}(\delta \mathbf{l}(\omega,x,n)^{-1})$, one has
    \begin{enumerate}[(i)]
        \item $\text{Lip}(\mathrm{d}_{\cdot}\varphi_{(\omega,x),n})\leq \mathbf{l}(\omega,x,n)$,
        \item $\text{Lip}(\varphi_{(\omega,x),n}-\mathrm{d}_0\varphi_{(\omega,x),n})<\delta$. 
    \end{enumerate}
\end{enumerate}
\end{proposition}
\begin{proof}
Recall the $\rho_b: M\to (0,1)$ be the regular radius for sufficient small $b$ as \Cref{26-1-22-1459}.  Let 
$$\mathbf{l}_3(\omega,x)=\frac{\max\{4\sup_{y\in B(x,1)}|\mathrm{d}_y\varphi_{\omega}|, 1\}}{\rho_b(x)\rho_b(\varphi_{\omega}(x))}.$$ 
By \Cref{26-7-26-1607}  we know that  $\exp_{\varphi_{\omega}(x)}^{-1}\circ \varphi_{\omega}\circ \exp_x : \mathcal{B}_x(\mathbf{l}_3(\omega,x)^{-1})\to\mathcal{B}_{\varphi_{\omega}x}(\rho_b(\varphi_{\omega}(x)))$ is well-defined.  By using the integrable condition,  we know that  
$$  \mathbf{l}_4(\omega,x):=2\mathbf{l}_3(\omega,x)\cdot\max\{\sup_{y\in B(x,1)}|\nabla_{y}\mathrm{d}\varphi_{\omega}|+\sup_{y\in B(x,1)}|\mathrm{d}_y\varphi_{\omega}|+\big(\sup_{y\in B(x,1)}|\mathrm{d}_y\varphi_{\omega}|\big)^2,1\}$$
 is a tempered function from  \cite[Propositions 4.3.3 and 4.1.3]{A}.  Therefore,  we can choose a $\tau/3$-slowly varying function $\tilde{\mathbf{l}}_4$ such that $0<\mathbf{l}_4\leq\tilde{\mathbf{l}}_4$. Next, we denote 
\begin{align}
\label{24-4-5-2115}
\mathbf{l}(\omega,x,n ):=4\max\{\mathbf{l}_2(\omega,x,n),\mathbf{l}_2(\omega,x,n+1)\}\tilde{\mathbf{l}}_{4}(\varPhi^n(\omega,x)),
\end{align}
which is a $\tau$-slowly varying function.  Assertions \eqref{23-12-18-2136-1} and \eqref{23-12-18-2136-2} follow from \Cref{23-12-18-21-37}. As for \eqref{23-12-18-2136-3},  by \Cref{26-1-22-1459} and \Cref{26-7-26-1607}, for any $w\in\hat{\mathbf{B}}:= \mathbf{B}(\delta \mathbf{l}(\omega,x,n)^{-1})$ one has that 
\begin{align*}
|\mathrm{d}_w\varphi_{(\omega,x),n}|\leq \sqrt{2}\frac{1+b}{1-b}\mathbf{l}_2(\omega,x,n)\sup_{y\in B(\varphi^n_{\omega}x,1)}|\mathrm{d}_{y}\varphi_{\theta^n\omega}|\leq \mathbf{l}(\omega,x,n),
\end{align*}
and 
\begin{align*}
    \text{Lip}(\varphi_{(\omega,x),n}-\mathrm{d}_0\varphi_{(\omega,x),n})&\leq \sup_{w\in \hat{\mathbf{B}}}\left\|\mathrm{d}_w\left(\varphi_{(\omega,x),n}-\mathrm{d}_0\varphi_{(\omega,x),n}\right)\right\|\\
    & = \sup_{w\in \hat{\mathbf{B}}}\left\|\mathrm{d}_w\varphi_{(\omega,x),n}-\mathrm{d}_0\varphi_{(\omega,x),n}\right\|\\
    &\leq \sup_{w\in \hat{\mathbf{B}}}\mathbf{l}(\omega,x,n)\|w\|<\delta. 
\end{align*}
This completes the proof of  \Cref{24-3-8-1417}.
\end{proof}

\subsection{Decomposition of neutralized Bowen ball}\label{s.010602}
The purpose of this subsection is  to decompose the neutralized Bowen ball   by unstable graphs under assuming that $(\varphi,\mu)$ has positive Lyapunov exponents.  In this subsection, we fix a small positive constant $\epsilon\ll 1$ and  a constant $\tau$ with 
\begin{align}
\label{26-2-23-1639}
0< \tau\leq \frac{\epsilon}{36D},
 \end{align}
 where $D$ is the dimension of Riemannian manifold $M$.  Through this subsection, we mainly consider  the neutralized Bowen ball in the $\tau$-Lyapunov chart  defined in \Cref{26-2-23-1409}. Firstly, we give a technical lemma.

\begin{lemma}
\label{25-2-11-2309}
For any  $(\omega,x)\in\mathbf{\Gamma}$,  for any sufficiently large $n\in\mathbb{N}$ (depending on $(\omega,x)$),  one has that 
\begin{align}
\label{24-8-20-2018}
\|\mathrm{d}_w\varphi_{(\omega,x),i}-\mathrm{d}_0\varphi_{(\omega,x),i}\|\leq e^{-n\epsilon/2}\leq\epsilon,
\end{align}
wherever  $w\in\mathbf{B}( e^{-n\epsilon})$ and $0\leq i\leq n$.
\end{lemma}
\begin{proof}
According to \Cref{24-3-8-1417} and \eqref{26-2-23-1639},  for any $0\leq i\leq n$ one has that 
\begin{align*}
\mathbf{l}(\omega,x, i)&\leq \mathbf{l}(\omega,x,0) e^{i\tau}\leq e^{n\epsilon/18D}
\end{align*}
for any sufficiently large $n$. By (c) in \Cref{24-3-8-1417}, we can obtain \Cref{25-2-11-2309}.
\end{proof}
Now, we would like to establish the pseudo-unstable  manifolds for  one iteration.  For any linear operator $L:\mathbf{R}^D\to\mathbf{R}^D$, we denote
\begin{align}
    \label{24-9-6-1010}
    |L|_{\mathrm{min}}=\inf_{v\neq 0}\frac{\|Lv\|}{\|v\|}.
\end{align}

\begin{lemma}
\label{25-12-17-1642}
  Fix a point  $(\omega,x)\in\mathbf{\Gamma}$  and a large enough $n\in\mathbb{N}$.  Let  $ r_1<e^{-n\epsilon}$ be  a small positive constant   and 
$$h_0: \mathbf{B}^{u}(r_1)\to \mathbf{B}^{cs}(r_1)$$ 
be a $C^1$  mapping satisfying $\|\mathrm{d}_{\cdot}h_0\|\leq 1$. Then for any $0\leq i\leq n-1$ the following properties hold. 
\begin{itemize}
\item The mapping 
$$g:=\pi^u\circ\varphi_{(\omega,x),i}\circ (\text{id}_u+h_0): \mathbf{B}^{u}(r_1)\to E^u$$
 is a diffeomorphism onto  with $|\mathrm{d}_{\cdot}g|_{\mathrm{min}}\geq e^{\chi/2}$ and $\mathbf{B}^{u}(r_1)\subset g(\mathbf{B}^{u}(r_1))$, where 
$\text{id}_u: E^u\to E^u$ is the identity map. 

\item There exists a $C^1$  mapping $h_1:\mathbf{B}^{u}(e^{3\epsilon}r_1)\to \mathbf{B}^{cs}(e^{3\epsilon}r_1)$ such that  
$$\|\mathrm{d}_{\cdot}h_1\|\leq1 \quad \text{ and }\quad \graph(h_1)\subset \varphi_{(\omega,x),i}\graph(h_0).$$
\end{itemize}
\end{lemma}
\begin{proof}
By (b)  in \Cref{24-3-8-1417}  and  \Cref{25-2-11-2309},  for any $w\in\mathcal{I}_0$ and $v\in E^u$, one has that
\begin{align}\label{e.051401}
\begin{split}
    &\|\mathrm{d}_{w}\big(\pi^u\circ \varphi_{(\omega,x),i}\circ(\text{id}_u+h_0)\big)v\|
    \\=&\|\pi^u\circ\mathrm{d}_{w+h_0(w)}\varphi_{(\omega,x),i}(v)+
    \pi^u\circ\mathrm{d}_{w+h_0(w)}\varphi_{(\omega,x),i}\circ\mathrm{d}_{w}h_0(v)\|
    \\\geq&\|\pi^u\circ\mathrm{d}_{w+h_0(w)}\varphi_{(\omega,x),i}(v)\|-\|
    \pi^u\circ\mathrm{d}_{w+h_0(w)}\varphi_{(\omega,x),i}\circ\mathrm{d}_{w}h_0(v)\|
    \\=&\|\pi^u\circ\mathrm{d}_{w+h_0(w)}\varphi_{(\omega,x),i}(v)\|-\|
    \pi^u\circ(\mathrm{d}_{w+h_0(w)}\varphi_{(\omega,x),i}-\mathrm{d}_0\varphi_{(\omega,x),i})\circ\mathrm{d}_{w}h_0(v)\|
    \\\geq&(e^{\chi}-\epsilon)\|v\|-\epsilon\|\mathrm{d}_{w}h_0v\|
    \\\geq& \big(e^{\chi}-2\epsilon\big)\|v\|
    \\\geq &e^{\chi/2}\|v\|,
\end{split}
\end{align}
where the last inequality holds  by choosing sufficiently small  $\epsilon$.  It follows that $g=\pi^u\circ\varphi_{(\omega,x),i}\circ (\text{id}_u+h_0)$ is  an injective mapping,  and $g$ is  a diffeomorphism onto its image by the inverse function theorem.  Next, we would like to verify its image cover $\mathbf{B}^u(r_1)$.  Note that 
\begin{align*}
\|\pi^u\circ \varphi_{(\omega,x),i}\circ h_0(0)\|&=\|\pi^u\circ \varphi_{(\omega,x),i}(h_0(0)-0)\|
\\&\leq \sup_{v\in \mathbf{B}^{cs}(r_1)}\big\|\pi^u\circ\mathrm{d}_v\varphi_{(\omega,x),i}|_{E^{cs}}\big\|\cdot\|h_0(0)\|
\\&\leq \sup_{v\in\mathbf{B}^{cs}(r_1)}\big\|\pi^u\circ\mathrm{d}_v\varphi_{(\omega,x),i}|_{E^{cs}}-\pi^u\circ\mathrm{d}_0\varphi_{(\omega,x),i}|_{E^{cs}}\big\| r_1
\\&\leq \epsilon r_1.
\end{align*}
It follows that 
\begin{align*}
g(\mathbf{B}^{u}(r_1))&\supset \mathbf{B}^u(g(0), e^{\chi/2}r_1)
\\&\supset \mathbf{B}^u(0, e^{\chi/2}r_1-\epsilon r_1)
\\&\supset \mathbf{B}^u( e^{3\epsilon}r_1)
\end{align*}
Therefore,  
$$g^{inv}:  \mathcal{I}_1:=\mathbf{B}^u( e^{3\epsilon}r_1)\to\mathbf{B}^u(r_1)\quad w\mapsto g^{-1}(w)$$ 
is a $C^1$ mapping  and $\|\mathrm{d}_{\cdot}g^{inv}\|\leq  e^{-\chi/2}$. Let
$$h_1: \mathcal{I}_1\to E^{cs},\quad u\mapsto \pi^{cs}\circ\varphi_{(\omega,x),i}\circ\big(\text{id}_u+h_0)\circ g^{inv}(u).$$ 
It is clear that $h_1$ is $C^1$ and $\varphi_{(\omega,x),i}\graph(h_0)\supset\graph(h_1)$. 

For any $w\in\mathcal{I}_0:=\mathbf{B}^{u}\big(r_1\big)$, similar to the derivation of \eqref{e.051401}, by (b) in \Cref{24-3-8-1417} and  \Cref{25-2-11-2309}, one has that
\begin{align*}
&\|\mathrm{d}_{w}\big(\pi^{cs}\circ\varphi_{(\omega,x),i}\circ\big(\text{id}_u+h_0)\circ g^{inv}\big)\|
\\\leq&\|\mathrm{d}_{g^{inv}(w)+h_0\circ g^{inv}(w)}\varphi_{(\omega,x),i}|_{E^{cs}}\circ\mathrm{d}_{g^{inv}(w)}h_0\circ\mathrm{d}_{w}g^{inv}\|+\|\pi^{cs}\circ\mathrm{d}_{g^{inv}(w)+h_0\circ g^{inv}(w)}\varphi_{(\omega,x),i}|_{E^{u}}\circ\mathrm{d}_{w}g^{inv}\|
\\\leq&e^{-\chi/2}\big(e^{\tau}+\epsilon+\epsilon\big)
\\<& e^{-\chi/3}<1,
\end{align*}
by choosing sufficiently small $\epsilon$.  
For any $w\in\mathbf{B}^u(r_1)$, one has that 
\begin{align*}
&\|\pi^{cs}\circ \varphi_{(\omega,x),i}\circ (w+h_0(w))\|
\\\leq&\|\pi^{cs}\circ \varphi_{(\omega,x),i}(w+h_0(w))-\pi^{cs}\circ \varphi_{(\omega,x),i}(w)\|+\|\pi^{cs}\circ \varphi_{(\omega,x),i}(w)-0\|
\\\leq &\sup_{v\in \mathbf{B}(r_1)}\|\pi^{cs}\circ\mathrm{d}_v\varphi_{(\omega,x),i}|_{E^{cs}}\|\cdot\|h_0(w)\|+\sup_{v\in \mathbf{B}(r_1)}\|\pi^{cs}\circ\mathrm{d}_v\varphi_{(\omega,x),i}|_{E^u}\|\cdot\|w\|
\\\leq &(e^{\tau}+\epsilon)\|h_0(w)\|+\epsilon\|w\|
\\\leq& (2\epsilon+e^{\tau})r_1
\\\tag{By \eqref{26-2-23-1639} }\leq &e^{3\epsilon}r_1.
\end{align*}
 This finishes the proof of  \Cref{25-12-17-1642}.
\end{proof}

Finally, we apply \Cref{25-12-17-1642} to decompose the neutralized Bowen ball in $\tau$-Lyapunov chart as  following 
\begin{align}
\label{26-2-23-1936}
\mathbf{B}(\omega,x, n, e^{-n\epsilon}):=\{w\in \mathbf{R}^d: \|\varphi^i_{(\omega,x)}w\|\leq e^{-n\epsilon} \text{ for } i=0,\dots,n-1\},
\end{align}
where $\varphi^i_{(\omega,x)}=\varphi_{(\omega,x),i-1}\circ\cdots\circ\varphi_{(\omega,x),0}$. The reason why we use the neutralized Bowen ball in \Cref{24-4-5-2040} is that this kind of Bowen  ball guarantees the existence of pseudo-unstable manifolds for the $n$-th iteration.
\begin{proposition}
\label{24-8-21-2303}
Given  a point   $(\omega,x)\in\mathbf{\Gamma}$, for any sufficiently large $n\in\mathbb{N}$ (maybe dependent on $(\omega,x)$)  and any $v\in\mathbf{B}^{cs}( e^{-4n\epsilon})$,    there  is an    open subset $\mathcal{I}_v$  of $\mathbf{B}^{u}(1)$  such that following mapping 
\begin{align}
g_v^{(n)}: \mathcal{I}_v\to \mathbf{B}^{u}(e^{-n\epsilon}),\quad w\mapsto \pi^u\circ\varphi^{n}_{(\omega,x)}\big(v+w\big)
\end{align}
is a diffeomorphism, and 
\begin{align}
\mathbf{B}(\omega,x, n, e^{-n\epsilon})\supset\bigsqcup_{v\in \mathbf{B}^{cs}( e^{-4n\epsilon})}\left(v+\mathcal{I}_v\right). 
\end{align}
\end{proposition}
\begin{proof}
Let 
$$h_v: \mathbf{B}^{u}(e^{-4n\epsilon})\to \mathbf{B}^{cs}(e^{-4n\epsilon}),\quad w\mapsto v.$$
 By applying \Cref{25-12-17-1642} for the mapping 
$h^0_v:=h_v$,  we infer the existence of a $C^1$ mapping $h_v^1: \mathbf{B}^u(e^{-4n\epsilon+3\epsilon})\to \mathbf{B}^{cs}(e^{-4n\epsilon+ 3\epsilon}),$ with the following properties 
\begin{itemize}
\item $\|\mathrm{d}_{\cdot} h^1_v\|\leq 1$;
\item $\text{graph}(h_v^{1})\subset\varphi_{(\omega,x),0}(\text{graph }h^0_v)$;
\item the mapping 
$$g_v^1:\mathcal{I}^{(1)}_v\to \mathbf{B}^u(r_1)\,\quad w\mapsto \pi^u\circ\varphi_{(\omega,x),0}(v+w)$$ is a diffeomorphism with $|\mathrm{d}_{\cdot}g_v^1|_{\mathrm{min}}\geq e^{\chi/2}$ , where $\mathcal{I}^{(1)}_v=(g_v^1)^{-1}\mathbf{B}^u(r_1)$. 
\end{itemize}
Repeating above arguments, we obtain two finite sequences of mappings 
\begin{align}
\label{24-10-29-1647}
    h_{v}^i: \mathcal{I}_v^{(i)}\to E^{cs}, \text{ and }g_{v}^i=\pi^u\circ \varphi^{}_{(\omega,x),i-1}\circ(\text{id}_u+h_{v}^{i-1}),
\end{align}
satisfying similar properties as above for any  $i=1,\dots,n$. In particular 
\begin{align}\label{e.051801}
    \|\mathrm{d}_{\cdot} h^i_v\|\leq 1, \text{ and } |\mathrm{d}_{\cdot}g_v^i|_{\mathrm{min}}\geq e^{\chi/2}.
\end{align}
The proof is then complete by setting $ g_v^{(n)} = g_{v}^{n}\circ g_{v}^{n-1} \circ \cdots \circ g_{v}^{1}$ and 
\[\mathcal{I}^{n}_v:=\mathcal{I}^{(n)}_v,\qquad \mathcal{I}^0_{v}= (g_{v}^{n})^{-1}\mathcal{I}^{(n)}_v.\]
This finishes the proof of  \Cref{24-8-21-2303}.
\end{proof}
In fact, by a direct calculation in the proof of \Cref{25-12-17-1642}, one has following remark.
\begin{remark}
\label{26-2-23-2100}
If the measure-preserving  RDS $(\varphi,\mu)$ has no positive Lyapunov exponents on $\mathbf{\Gamma}$, then  for any at   $(\omega,x)\in \mathbf{\Gamma}$, for any sufficiently large $n\in\mathbb{N}$ one has that 
$$\mathbf{B}(e^{-4n\epsilon})\subset \mathbf{B}(\omega,x, n, e^{-n\epsilon}).$$
\end{remark}


\subsection{Volume  estimations  of neutralized Bowen ball}
\label{24-10-29-2033}
In this subsection, we would like to  estimate the volume of  neutralized Bowen ball defined in \eqref{26-2-23-1936} with respect to the Lebesgue measure. Throughout this section, $C$ denotes a generic constant that is independent of $n$ and $\epsilon$, but may vary from line to line. Continuing with the notation from \Cref{24-8-21-2303}, we aim to prove the following 
\begin{proposition}
\label{24-7-19-2157}
  Fix a $(\omega,x)\in\mathbf{\Gamma}$. There  exists a constant $C>0$ such that for sufficiently large   $n\in\mathbb{N}$  one has that 
\begin{align}
e^{-Cn\epsilon}\leq \textbf{vol}\Big( \mathbf{B}(\omega,x, n, e^{-n\epsilon})\Big)\cdot\det(\mathrm{d}_0\varphi^{n}_{(\omega,x)}|_{E^u}),
\end{align}
where  $\textbf{vol}$ is the Lebesgue measure on $\mathbf{R}^D$.  
\end{proposition}

Before  proving \Cref{24-7-19-2157}, we give two  technical lemmas  to help us to  relate  the volume expansion rate along  the expansion rate of $g_v^{(n)}$ to the unstable subspace. 
  \begin{lemma}
    \label{24-4-6-1027}
    For any $D\in\mathbb{N}$ and  $C>1$,  if $A, B:\mathbf{R}^D\to\mathbf{R}^D$ are invertible bounded linear operators and for which
    \begin{align*}
    |A|,|B|,|A^{-1}|, |B^{-1}|\leq C,
    \end{align*}
    then 
    $$\left|\log\frac{\det(A)}{\det(B)}\right|\leq DC^{2D-1}|A-B|.$$
    \end{lemma}
    \begin{proof}
    Let $e_1,\dots, e_D$ be an orthonormal basis of $\mathbf{R}^D$.  Then one has that 
    \begin{align*}
    &|\det(A)-\det(B)|
    \\=&\big||A(e_1\wedge e_2\cdots\wedge e_D)|-|B(e_1\wedge e_2\cdots\wedge e_D)|
    \big|
    \\\leq& |A(e_1\wedge e_2\cdots\wedge e_k)-B(e_1\wedge e_2\cdots\wedge e_D)|
    \\\leq& |A-B|\cdot|A(e_2\wedge e_2\cdots\wedge e_D)|+|B|\cdot |A(e_2\wedge e_3\cdots\wedge e_D)-B(e_2\wedge e_3\cdots\wedge e_D)|
    \\\leq& |A|^{D-1}\cdot |A-B|+|A(e_2\wedge e_3\cdots\wedge e_D)-B(e_2\wedge e_3\cdots\wedge e_D)|\cdot|B|
    \\\leq&\sum_{i=0}^{D-1}|A|^{D-1-i}\cdot|A-B|\cdot|B|^i
    \\\leq& DC^{D-1}|A-B|.
    \end{align*}
    By the mean value theorem and 
    $$\min\{\det(A), \det(B)\}\geq \min\{|A^{-1}|^{-D}, |B^{-1}|^{-D}\}\geq C^{-D}, $$  one obtains 
    \begin{align*}
    \big|\log\det(A)-\log\det(B)\big|&\leq C^D|\det(A)-\det(B)|
    \\&\leq DC^{2D-1}|A-B|.
    \end{align*}
    \end{proof}

\begin{lemma}
\label{24-7-19-2244}
Fix a point  $(\omega,x)\in\mathbf{\Gamma}$. There  exists a constant $C>1$ such that for sufficiently large   $n\in\mathbb{N}$, if  $v\in\mathbf{B}^{cs}( e^{-4n\epsilon})$  and $w\in\mathcal{I}_{v}$, then 
\begin{align}
& C^{-1}\leq\frac{\det(\mathrm{d}_0\varphi^{n}_{(\omega,x)}|_{E^u})}{\det\big(\mathrm{d}_{w}g^{(n)}_v\big)}\leq C.
\end{align}
\end{lemma}
\begin{proof}
Recall that  for each $i=1,\dots,n$
$$h_{v}^{i-1}: \mathcal{I}_v^{(i)}\to E^{cs} \text{ and } g_{v}^i=\pi^u\circ \varphi_{(\omega,x),i-1}\circ(\text{id}_u+h_{v}^{i-1})$$
are defined as in \eqref{24-10-29-1647} in the proof of  \Cref{24-8-21-2303}. For $w\in\mathcal{I}_v$ we have 
\[w_{i-1}+h_v^{i-1}(w_{i-1})\in\mathbf{B}(e^{-n\epsilon}), \quad i=1,\dots, n,\]
by denoting $w_0:=w$,  and   
\[ w_i=g_v^i\circ\cdots\circ g_v^{1}(w), \quad i\geq 1.\]
It also follows from the chain rule that 
$$\mathrm{d}_{w}g_v^{(n)}=\mathrm{d}_{w_{n-1}}g_v^n\circ\cdots \circ\mathrm{d}_{w_0}g_v^1.$$

For $i=1,\dots, n$, we infer from \Cref{24-3-8-1417}, and \eqref{26-2-23-1639}  that 
\begin{align}\label{24-4-14-1136}
\begin{split}
    &\|\mathrm{d}_{w_{i-1}}g_v^i\|=\|\mathrm{d}_{w_{i-1}}\big(\pi^u\circ \varphi^{}_{\varPhi^{i-1}(\omega,x)}\circ(\text{id}_u+h_{v}^{i-1})\big)\|\leq e^{n\epsilon/9D},\\
    &\|\mathrm{d}_0\varphi^{}_{\varPhi^{i-1}(\omega,x)}|_{E^u}\|\leq e^{n\epsilon/9D}.
\end{split}
\end{align}
by letting  $n$ big enough.
Besides,  by \eqref{e.051801} and  \Cref{24-3-8-1417}, one has
\begin{align}
\label{24-4-14-1137}
|\mathrm{d}_{w_{i-1}}g_v^i|_{\mathrm{min}}\geq e^{\chi/2}\quad \text{and}\quad |\mathrm{d}_0\varphi^{}_{\varPhi^{i-1}(\omega,x)}|_{E^u}|_{\mathrm{min}}\geq e^{\chi/2}.
\end{align}
 By the invariance properties in item (b) of Proposition \ref{24-3-8-1417} and \eqref{24-8-20-2018}, one also has
\begin{align}\label{e.051802}
    \begin{split}
        &\left\|\mathrm{d}_{w_{i-1}}g_v^i - \mathrm{d}_0\varphi^{}_{\varPhi^{i-1}(\omega,x)}\mid_{E^u}\right\|\\
        = &\sup_{u\in \mathbf{B}^u(1)}\left\|\pi^u\circ\mathrm{d}_{\widetilde{w}_{i-1}}\varphi^{}_{\varPhi^{i-1}(\omega,x)}(u+\mathrm{d}_{w_{i-1}}h_{v}^{i-1}u) - \pi^u\circ\mathrm{d}_0\varphi^{}_{\varPhi^{i-1}(\omega,x)}(u+ \mathrm{d}_{w_{i-1}}h_{v}^{i-1}u)\right\|\\
         \leq &e^{-n\epsilon/2},
    \end{split}
\end{align}
where $\widetilde{w}_{i-1}=w_{i-1}+h_v^{i-1}(w_{i-1})$. 

\medskip

Therefore, by  \Cref{24-4-6-1027},  \eqref{24-4-14-1136}, \eqref{24-4-14-1137} and \eqref{e.051802},  we have
\begin{align}
\notag&|\log\det(\mathrm{d}_{w}g_v^{(n)})-\log\det(\mathrm{d}_0\varphi^{n}_{(\omega,x)}|_{E^u})|
\\\notag\leq&\sum_{i=0}^{n-1}|\log\det(\mathrm{d}_{w_{i}}g_v^{i+1})-\log\det(\mathrm{d}_0\varphi^{}_{\varPhi^{i}(\omega,x)}|_{E^u})|
\\\notag\leq&\sum_{i=0}^{n-1}D  e^{n\epsilon/3}\left\|\mathrm{d}_{w_{i}}g_v^{i+1}-\mathrm{d}_0\varphi^{}_{\varPhi^{i}(\omega,x)}|_{E^u}\right\|
\\\notag\leq&\sum_{i=0}^{n-1}D   e^{n\epsilon/3} e^{-n\epsilon/2}
\\\notag\leq &  C,
\end{align}
by choosing $n$ large enough. This completes the proof of \Cref{24-7-19-2244}.
\end{proof}

Let's proceed to   prove \Cref{24-7-19-2157}. 
\begin{proof}
For any $v\in\mathbf{B}^{cs}( e^{-4n\epsilon})$, one has that 
\begin{align*}
&\text{vol}_u(\mathcal{I}_{v})\cdot\det(\mathrm{d}_0\varphi^{n}_{(\omega,x)}|_{E^u})
\\=&\int_{E^u}1_{\mathcal{I}_{v}}(w)\mathrm{d}\text{vol}_u(w)\cdot\det(\mathrm{d}_0\varphi^{n}_{(\omega,x)}|_{E^u})
\\\tag{by \Cref{24-7-19-2244}}\geq & C^{-1}\int_{E^u}1_{\mathcal{I}_{v}}(w)\cdot\det(\mathrm{d}_{w}g_v^{(n)})\mathrm{d}\text{vol}_u(w)
\\\tag{by \Cref{24-8-21-2303}}\geq&C^{-1}\int_{\mathbf{B}^u(e^{-n\epsilon})}\mathrm{d}\text{vol}_u(w)
\\\geq &e^{-Cn\epsilon},
\end{align*}
where $\text{vol}_u$ is the  Lebesgue  measure on $E^u$ and $n$ is large enough.  Together with \Cref{24-8-21-2303} and Fubini theorem, we can obtain  \Cref{24-7-19-2157}. 
\end{proof}

\subsection{Proof of \Cref{24-4-5-2040}}
\label{26-2-23-1418}
In this subsection, we give the final  proof of  \Cref{24-4-5-2040}. Denote $\text{vol}$ as the Lebesgue measure on $M$ induced by its Riemannian structure.  Due to  \Cref{24-7-19-2157} and \eqref{25-11-12-1407}, it is not hard to see that for $\mu$-a.s. $(\omega,x)\in\Omega\times M$ one has that 
\begin{align}
\label{25-2-12-2135}
\lim_{\epsilon\to 0}\limsup_{n\to\infty}-\frac{\log\text{vol}(B(\omega,x, n, e^{-n\epsilon}))}{n}\leq\lim_{n\to\infty}\frac{\log\det(\mathrm{d}_x\varphi^n_{\omega}|_{E_0^u(\omega,x)})}{n}= \lambda^+_{\sum}(\omega,x).
\end{align}
It is worth mentioning that  \eqref{25-2-12-2135} still holds when $(\varphi,\mu)$ has no positive Lyapunov exponents by  \Cref{26-2-23-2100}.  
\begin{proof}[Proof of  \Cref{24-4-5-2040}]
Suppose that \Cref{24-4-5-2040}  doesn't hold.  Then there exists a positive $\mu$-measure subset $A$  such that for any $(\omega,x)\in A$ such that 
\begin{align*}
\lim_{\epsilon\to 0}\limsup_{n\to\infty}-\frac{\log\mu_{\omega}(B(\omega,x, n, e^{-n\epsilon}))}{n}> \lambda^+_{\sum}(\omega,x).
\end{align*} 
By \eqref{25-2-12-2135}, we also have that 
\begin{align}
\lim_{\epsilon\to 0}\limsup_{n\to\infty}-\frac{\log\text{vol}(B(\omega,x, n, e^{-n\epsilon}/3))}{n}\leq \lambda^+_{\sum}(\omega,x).
\end{align}
By Lusin's lemma, there exists a  small constant $\delta_0$, a positive constant $\lambda$ and a positive $\mu$-measurable subset $A_1\subset A$ such that  for any $(\omega,x)\in A_1$ one has that 

\begin{align*}
\lim_{\epsilon\to 0}\limsup_{n\to\infty}-\frac{\log\mu_{\omega}(B(\omega,x, n, e^{-n\epsilon}))}{n}&> \lambda+2\delta_0
\\>\lambda-\delta_0&\geq \lim_{\epsilon\to 0}\limsup_{n\to\infty}-\frac{\log\text{vol}(B(\omega,x, n, e^{-n\epsilon}/3))}{n} .
\end{align*} 
By Egorov's theorem, we can choose small constant $\epsilon>0$  and  a positive $\mu$-measurable subset $A_2\subset A_1$ such that for any $(\omega,x)\in A_2$ one has that 
\begin{align*}
\limsup_{n\to\infty}-\frac{\log\mu_{\omega}(B(\omega,x, n, e^{-n\epsilon}))}{n}> \lambda+\delta_0\geq\lambda\geq \limsup_{n\to\infty}-\frac{\log\text{vol}(B(\omega,x, n, e^{-n\epsilon}/3))}{n}.
\end{align*} 
By using Egorov's theorem again,  we can choose  a compact  positive $\mu$-measurable subset $A_3\subset A_2$ and a large number $N_0$  such that  any $N\geq N_0$ and $(\omega,x)\in A_3$ one has that  
$$\text{vol}(B(\omega,x, N, e^{-N\epsilon}/3))\geq e^{-N(\lambda+\delta_0/4)}.$$
Fixing $N$ sufficiently large, denote for any $n\geq N$ 
$$A_{3, n}:=\{(\omega,x)\in A_3 :  n  \text{ is the smallest integer s.t. } \mu_{\omega}(B(\omega,x,n, e^{-n\epsilon}))\leq e^{-n(\lambda+3\delta_0/4)}\}.$$
For any  $\omega\in\pi_{\Omega}(A_{3,n})$, we can choose a  sequence of points of $\{x_i: i=1,\dots, s_n\}$ such that  $\{B(\omega, x_i, n, e^{-n\epsilon}): i=1,\dots, s_n\}$ spans the subset $A_{3,n}(\omega)$ and $\{B(\omega,x_i,n, e^{-n\epsilon}/3): i=1,\dots, s_n\}$ is a separated set  of  $A_{3,n}(\omega)$. It is clear that  $s_n\leq C_K e ^{n(\lambda+\delta_0/4)}$, where $C_K$ is a constant depending on the volume 
$$K:=\{y\in M: d\big(y, \pi_M(A_3)\big)\leq 1\}.$$  Since the balls $B(\omega,x_i,n,e^{-n\epsilon})$ cover $A_{3,n}(\omega)$ and each center $x_i$ belongs to $A_{3,n}(\omega)$, we have
\begin{align*}
\mu_\omega(A_{3,n}(\omega))&\leq \sum_{i=1}^{s_n}\mu_\omega(B(\omega,x_i,n,e^{-n\epsilon}))
\\&\leq s_n e^{-n(\lambda+3\delta_0/4)}
\\&\leq C_Ke^{n(\lambda+\delta_0/4)}e^{-n(\lambda+3\delta_0/4)}=C_Ke^{-n\delta_0/2}.
\end{align*}
Therefore, 
$$\mu_{\omega}(A_3(\omega))\leq\sum_{n\geq N}\mu_{\omega}(A_{3,n}(\omega))\leq \frac{C_Ke^{-N\delta_0/2}}{1-e^{-\delta_0/2}},$$
which implies that $\mu(A_3(\omega))=0$ when $N$ tends to infty.  This leads to a contradiction.  This completes the proof.
\end{proof}

\section{Counterexample of Brin-Katok local entropy formula}
\label{26-5-26-0943}
In this section, we construct a counterexample showing that the integrability
condition \eqref{2026-2-15-1654} in \Cref{26-2-13-1805} cannot in general be
omitted. More precisely, we construct a smooth diffeomorphism on a complete
non-compact Riemannian manifold and an ergodic invariant probability measure
for which the measure-theoretic entropy is finite, whereas the upper
Brin--Katok local entropy is infinite almost everywhere. 

\begin{thm}\label{thm:main}
There exists a connected, complete, non-compact smooth Riemannian manifold $M$, a smooth diffeomorphism $f:M\to M$, and an $f$-ergodic  measure $\mu$ such that

$$ h_\mu(f)<\infty, \text{ but } \,  \overline h_\mu(f,x)=+\infty,$$
for $\mu$-a.s. $x\in M$.

\end{thm}

\subsection{Example of Construction}
\label{26-7-8-1617}
Denote  $\mathbb{T}^2:=\mathbb{R}^2/\mathbb{Z}^2,$ and 
$$
  A:=\begin{pmatrix}2&1\\1&1\end{pmatrix}.
$$
Let $F:\mathbb{T}^2\to\mathbb{T}^2$ be the linear torus automorphism induced by $A$.
Let $\text{g}_0$ be the standard flat Riemannian metric on $\mathbb{T}^2$, $d_0$ the  distance induced by $\text{g}_0$, and $m$ the normalized Lebesgue area measure.  It is well known that 
\begin{enumerate}[(a)]
\item $m$ is an ergodic-measure of $F$ and  the measure-theoretic $h_{m}(F)=\log\frac{3+\sqrt{5}}{2}.$
\item  $(F,m)$ is exponential decay of correlations for H\"older observables. 
\end{enumerate}

Define
\[
  M:=\mathbb{T}^2\setminus\{0\},
  \qquad
  f:=F|_M:M\to M,
  \qquad
  \mu:=m|_M.
\]
Since  $0$ is a fixed point and $m(0)=0$,  $f$ is a smooth diffeomorphism of $M$, $\mu$ is an $f$-ergodic  measure, and
$$ h_\mu(f)=h_m(F)=\log\frac{3+\sqrt{5}}{2}.$$

We now start to define a new   Riemannian metric on $M$.  Fix a small constant $\eta>0$ and  a distance function $\mathbf{r}(x):=d_0(x,0)$.  Then we can 
 choose a smooth cut-off function $\chi:M\to[0,1]$ such that
$$
  \chi(x)=1 \quad (0<\mathbf{r}(x)<\eta/3),
  \qquad
  \chi(x)=0 \quad (\mathbf{r}(x)>2\eta/3).
$$
It follows that 
  $$\mathbf{a}(x):=\exp\!\left(\chi(x)e^{\mathbf{r}(x)^{-2}}\right),$$
is  a smooth positive function on $M$.  
\begin{lemma}\label{lem:complete}
Let $\text{g}:=\mathbf{a}^2\text{g}_0$. Then $(M,\text{g})$ is a non-compact complete Riemannian manifold.
\end{lemma}

\begin{proof}
Obviously, $M=\mathbb{T}^2\setminus\{0\}$ is non-compact.  We only need to  show that $\text{g}$ is complete. Denote $d_{\text{g}}$ as the Riemannian metric induced by $\text{g}$.
For  any $d_{\text{g}}$-Cauchy sequence $(x_n)_{n\in\mathbb{N}}$, one has that  it is also a $d_0$-Cauchy sequence since $\mathbf{a}\geq1$. So it converges to some point $q$ on the  torus $\mathbb{T}^2$. If $q\ne 0$, then on a relatively compact neighbourhood of $q$, $\mathbf{a}$ is bounded above and below, which implies that  $d_{\text{g}}$ is locally bi-Lipschitz equivalent to $d_0$. Thus $(x_n)_{n\in\mathbb{N}}$ converges to  $q$ under the metric $d_{\text{g}}$.

If $q=0$,  we claim that $\lim_{n\to\infty}d_{\text{g}}(x_n,x_1)=+\infty$,  contradicting the property of  Cauchy sequence.  Without loss of generality,  we assume that  the sequence $(x_n)_{n\in\mathbb{N}}$ lies a small neighborhood of $0$. Let $\gamma:[0,1]\to M$ be a  piecewise $C^1$ cure   on $\mathbb{T}^2$ which connects  $x_n$ and $x_1$.
Then one has that 
\begin{align*}
  \ell_{\text{g}}(\gamma)&=\int_0^1 \mathbf{a}(\gamma(t))\,|\dot\gamma(t)|_{\text{g}_0}\,\mathrm{d}t
  \\&=\int_0^1 \mathbf{a}(\gamma(t))\,|\dot\gamma(t)|_{\text{g}_0}\cdot|\nabla \mathbf{r}(\gamma(t))|_{\text{g}_0}\,\mathrm{d}t
 \\\tag{by $\tau=\mathbf{r}\circ\gamma$}& \geq \int_{\mathbf{r}(x_n)}^{\mathbf{r}(x_1)} \exp(e^{\tau^{-2}}) \mathrm{d}\tau,
\end{align*}
which implies the claim. By the Hopf--Rinow theorem, the Riemannian manifold $(M,\text{g})$ is complete.
\end{proof}

Finally, we will verify that $(f,\mu)$ doesn't satisfy condition \eqref{2026-2-15-1654}. Let  $E^u$ and $E^s$ be the unstable and stable space of $A$, respectively,  $\lambda$ be the eigenvalue of unstable space. In a sufficiently small coordinate neighborhood of $0$, write
$x=x^u+x^s$ according to the splitting $\mathbb R^2=E^u\oplus E^s$, and set
\[
\mathcal C_\lambda
:=
\{x\in M: |x^u|_{\text{g}_0}\le \lambda^{-2}|x^s|_{\text{g}_0}\}.
\]
 For $x \in \mathcal{C}_{\lambda}$ with $0 < \mathbf{r}(x) < \eta/3$,  it is obvious that $\mathbf{r}(f(x))\leq \sqrt{2}\lambda^{-1}\mathbf{r}(x)$.
Letting  $v\in E^s$  with $|v|_{\text{g}_0}=1$,  then one has that 
\begin{align*}
\log\|\mathrm{d}_xf\|_{\text{g}} &\geq\log\frac{\|\mathrm{d}_xf(v)\|_{\text{g}}}{\|v\|_{\text{g}}}
\\&= \log\Big(\frac{\mathbf{a}(fx)}{\mathbf{a}(x)} \cdot \|\mathrm{d}_xf(v)\|_{\text{g}_0}\Big)
\\&=e^{1/\mathbf{r}(f(x))^2}-e^{1/\mathbf{r}(x)^2}-\log\lambda
\\&\geq \frac{1}{2}e^{\lambda^2/2\mathbf{r}(x)^2}
\end{align*}
for $\mathbf{r}(x)$ sufficiently small.   Fixed small constant $\epsilon>0$, denote
$$
\mathcal{C}_{\lambda}^\epsilon = \{ x \in \mathcal{C}_{\lambda}: 0 < \mathbf{r}(x) < \epsilon \}.
$$
Using the  polar coordinates transformation, there exists a constant $C_{\epsilon}>0$ such that 
\[
\begin{aligned}
\int_{\mathcal{C}_{\lambda}^\epsilon}  \log^+\|\mathrm{d}_xf\|_{\text{g}}\,\mathrm{d}\mu(x)&\geq C_{\epsilon} \int_{0}^{\epsilon} e^{\lambda^2 / 2\tau^2} \tau\mathrm{d}\tau
\\&= \frac{C_{\epsilon}}{4}\int_0^{2\epsilon^2}e^{\lambda^2/s}\mathrm{d}s=+\infty.
\end{aligned}
\]
It follows  for any $b>0$ one has 
$$\int_{M}\log^+\sup_{y\in B_{d_{\text{g}}}(x,b)}\|\mathrm{d}_{y}f\|_{\mathrm{g}}\,\mathrm{d}\mu(x)\geq \int_{\mathcal{C}_{\lambda}^\epsilon}  \log^+\|\mathrm{d}_xf\|_{\text{g}}\,\mathrm{d}\mu(x)=+\infty.$$

\subsection{Two auxiliary lemmas}

In this subsection, we give  two auxiliary lemmas for the dynamical system $(M,f,\mu)$ defined  as \Cref{26-7-8-1617}. The first lemma can be regard as a kind of dynamical Borel--Cantelli lemma (for example, see \cite{MR3035294}). 

\begin{lemma}
\label{bc}
For any $k\geq 2$, let  $\rho_k:=(2\log k)^{-1/2}$ and $ E_k:=B_{d_0}(0,\rho_k).$ 
Then for $\mu$-a.s.  $x\in M$, there exist infinitely many $k$ such that  $f^k x\in E_k.$
\end{lemma}

\begin{proof}
Consider the measure-preserving dynamical system $(\mathbb{T}^2,F, m)$.  For  sufficiently large $k$, one has that
$$ m(E_k)=\pi\rho_k^2=\frac{\pi}{2\log k}\geq \frac{1}{\sqrt{k}}.$$
Since  $(\mathbb{T}^2,F, m)$ is exponential decay of correlations for H\"older observables, in particular for Lipschitz observables. Thus we can apply  \cite[Theorem 3.1]{MR3035294}  to the sequence of target balls $(E_k)_{k\geq 2}$ and yields
\[
 \lim_{n\to\infty} \frac{\sum_{k=1}^n \mathbf{1}_{E_k}(F^kx)}{\sum_{k=1}^n m(E_k)}=1
\]
for $m$-a.s. $x\in\mathbb{T}^2$.  It follows that $F^kx\in E_k$ occurs infinitely often. Removing the point $0$ does not affect the conclusion on a null set, so for $\mu$-a.s. $x\in M$, $f^kx\in E_k$ infinitely often.
\end{proof}

The second lemma is an estimate of  volume  of ordinary Riemannian balls under the new Riemanninan metric.

\begin{lemma}
\label{26-5-25-2035}
For any  $\epsilon>0$, there exists constant $r_\epsilon>0$ such that for any $z\in \{x\in M: \mathbf{r}(x)<r_{\epsilon}\}$, one has that 
 $$\mu\bigl(B_{d_{\text{g}}}(z,\epsilon)\bigr)\leq  16\pi\epsilon^2\mathbf{a}(z)^{-2}.$$
\end{lemma}

\begin{proof}
Given $\epsilon>0$, claim that  there is a  $r_\epsilon>0$ such that for any $z\in \{x\in M: \mathbf{r}(x)<r_{\epsilon}\}$ we can ensure that for any $w\in B_{d_0}(z, \frac{4\epsilon}{\mathbf{a}(z)})$, one has that
$$\mathbf{a}(w)\geq \frac{\mathbf{a}(z)}{2}.$$
Choosing $r_\epsilon>0$ smaller if necessary, we may assume that $\frac{8\epsilon}{\mathbf{a}(z)}<\mathbf{r}(z)$
whenever $\mathbf r(z)<r_\epsilon$.  Hence that  for any $w\in B_{d_0}(z, \frac{4\epsilon}{\mathbf{a}(z)})$, we have
$$
  \frac{\mathbf{r}(z)}{2}<\mathbf{r}(w)<\frac{3\mathbf{r}(z)}{2}.
$$
Let  $l(s):=e^{s^{-2}}.$  By  the mean value theorem, there exists a constant $C_1>0$ such that
$$
\begin{aligned}
  |\log\mathbf{a}(w)-\log\mathbf{a}(z)|
  &\le \sup_{\mathbf{r}(z)/2\leq s\leq 3\mathbf{r}(z)/2} |l'(s)|\,\frac{8\epsilon}{\mathbf{a}(z)} \\
  &\le C_1 \mathbf{r}(z)^{-3}e^{4/\mathbf{r}(z)^2}\cdot 8\epsilon \exp (-e^{1/\mathbf{r}(z)^2}).
\end{aligned}
$$
Since $\mathbf{r}(z)^{-3}e^{4/\mathbf{r}(z)^2} \exp (-e^{1/\mathbf{r}(z)^2})\to 0$ as $\mathbf{r}(z)\to 0$,    we can choose  small constant  $r_\epsilon>0$ such that the claim holds. 

Now,  fix $z\in \{x \in M: \mathbf{r}(x)<r_{\epsilon}\}$  and  $y\notin B_{d_0}(z,\frac{4\epsilon}{\mathbf{a}(z)})$. Take any piecewise $C^1$ curve $\gamma: [0,1]\to M$ joining $z$ and $y$. Let $\tau$ be the first time $\gamma$ leaves $B_{d_0}(z,\frac{4\epsilon}{\mathbf{a}(z)})$. Then, 
\[
  \ell_{\text{g}}(\gamma) \geq \frac{\mathbf{a}(z)}{2} \,\ell_{\text{g}_0}(\gamma|_{[0,\tau]})
  \ge \frac{\mathbf{a}(z)}{2} \frac{4\epsilon}{\mathbf{a}(z)}
  =2\epsilon.
\]
It follows that  $d_{\text{g}}(z,y)>\epsilon$. Particularly, one has that  $B_{d_\text{g}}(z,\epsilon)\subset B_{d_0}(z, \frac{4\epsilon}{\mathbf{a}(z)}).$
Hence that 
$$
  \mu(B_{d_{\text{g}}}(z,\epsilon))\leq 16\pi\epsilon^2\mathbf{a}(z)^{-2}.
$$
\end{proof}

\subsection{Estimation of upper Brin--Katok local entropy}

In this subsection,  we will verify that the upper Brin-Katok local entropy of $(M,f,\mu)$ is infinite under the new Riemannian metric. Fix any $\epsilon>0$.  Let  $\rho_k:=(2\log k)^{-1/2}$ for any $k\geq 2$.  By Lemma \ref{bc}, for $\mu$-almost every $x\in M$ there exist infinitely many $k$ such that
\[
  f^kx\in B_{d_0}(0, \rho_k).
\]
Fix such a point  $x$ and take an infinite subsequence $k\to\infty$ with this property.  Setting $z_k:=f^kx$, then
\begin{align}
\label{26-5-26-0953}
\mathbf{r}(z_k)<\rho_k=(2\log k)^{-1/2}.
\end{align}
For sufficiently large $k$,  it is clear that 
\[
  \log \mathbf{a}(z_k)=e^{1/\mathbf{r}(z_k)^2}>e^{2\log k}=k^2.
\]
Define the Bowen ball  as 
\[
  B(x, k+1,\epsilon)
  :=\{y\in M:\ d_{\text{g}}(f^jx,f^j y)<\epsilon\text{ for  } j=0,1,\cdots, k\}.
\]
Since  $\mu$ is $f$-invariant  and $f$ is bijective, one has that
\[
\begin{aligned}
  \mu\bigl(B(x,k+1,\epsilon)\bigr)&=\mu(f^k\bigl(B(x, k+1, \epsilon)\bigr))
\\&\leq \mu\bigl(B_{d_{\text{g}}}(f^kx,\epsilon)\bigr).
\end{aligned}
\]
By \eqref{26-5-26-0953} and \Cref{26-5-25-2035}, there exists  infinitely many sufficiently large $k$ such that 
\[
  \mu\bigl(B(x,k+1, \epsilon)\bigr)
  \le C_\epsilon \mathbf{a}(f^kx)^{-2},
\]
where $C_\epsilon:=16\pi\epsilon^2$.
By a direct  calculation, one has that 
\[
\begin{aligned}
  -\frac1{k+1}\log \mu\bigl(B(x,k+1, \epsilon)\bigr)
  &\ge \frac{2\log \mathbf{a}(f^kx)-\log C_\epsilon}{k+1}\\
  &\ge \frac{2k^2-\log C_\epsilon}{k+1}.
\end{aligned}
\]
 Consequently,  for $\mu$-a.s. $x\in M$ one has that 
\[
 \overline{h}_{\mu}(f,x)=\lim_{\epsilon\to0} \limsup_{n\to\infty} -\frac1n\log\mu(B(x,n, \epsilon))=+\infty.
\]

 \appendix 

\section{Geometric estimation}

Throughout this appendix,  we always assume that $M$ is a complete Riemannian manifold, all operator and tensor norms are induced
by the Riemannian metric. Each tangent space $T_xM$ is equipped with
its canonical flat connection, and pullback bundles over $T_xM$ are
equipped with the connection induced by the Levi--Civita connection
of $M$. 

\begin{lemma}
\label{26-7-26-1608}
 Given a small positive constant $b<1$,  there exists a measurable  regular radius  $\rho_{b}: M\to (0,1)$ such that for any $x\in M$
\begin{itemize}
\item $\exp_x: \mathcal{B}_x(\rho_b(x))\to B(x,1)$ is a diffeomorphism onto its image, where   $\mathcal{B}_x(\rho_b(x))$ is the open ball in  $T_xM$ at the  origin with  radius $\rho_b(x)$;
\item  $|\mathrm{d}_{.}\exp_x|\in (1-b, 1+b)$ and $|\nabla_{.}(\mathrm{d}\exp_x)|\leq b $ on $\mathcal{B}_x(\rho_b(x))$;
\end{itemize}
\end{lemma}
\begin{proof}
Consider the smooth map
$$
\mathcal{E}:TM\to M\times M,
\qquad
\mathcal{E}(y,v)=\bigl(y,\exp_y(v)\bigr).
$$
Using a local trivialization of $TM$ near a point $x\in M$, the differential of
$\mathcal{E}$ at $(x,0)$ is
\[
\mathrm{d}_{(x,0)}\mathcal{E}(\xi,\eta)=(\xi,\xi+\eta),
\]
which is  an isomorphism. By the inverse function theorem,
there exists an open neighborhood $\mathcal U_x$ of $(x,0)$ in
$TM$ such that $\mathcal{E}|_{\mathcal {U}_x}$ is a diffeomorphism onto its
image. Since
$\mathrm{d}_0\exp_y=\operatorname{id}_{T_yM}$ and $\nabla_0(\mathrm{d}\exp_y)=0$
for every $y\in M$, after shrinking $\mathcal U_x$ if necessary, we
may assume that
\[
(1-b)|\xi|<|\mathrm{d}_v\exp_y(\xi)|<(1+b)|\xi| \quad\text{and}\quad|\nabla_v(\mathrm{d}\exp_y)|\leq b,
\]
    whenever $(y,v)\in\mathcal W_x$ and $\xi\in T_yM$. 
    
We can choose an open neighborhood $U_x$ of $x$ and a number
$r_x\in(0,1)$ such that
\[
\bigl\{(y,v)\in TM:y\in U_x,\ |v|<r_x\bigr\}
\subset \mathcal W_x.
\]
For each $y\in U_x$, the restriction $\exp_y:\mathcal B_y(r_x)\to M$
is injective. Indeed, if $v,w\in\mathcal B_y(r_x)$ and
$\exp_y(v)=\exp_y(w)$, then
$$\mathcal{E}(y,v)=\mathcal{E}(y,w),$$
and the injectivity of $\mathcal{E}|_{\mathcal W_x}$ implies $v=w$. Thus $\exp_y|_{\mathcal{B}_y(r_x)}$ is an injective local
diffeomorphism, and hence a diffeomorphism onto its image. Furthermore, for every $v\in\mathcal B_y(r_x)$, the radial geodesic
\[
\gamma(t)=\exp_y(tv),\qquad 0\leq t\leq1,
\]
has length $|v|$. Consequently,
\[
d\bigl(y,\exp_y(v)\bigr)
\leq \ell(\gamma)
=|v|
<r_x<1.
\]
It follows that $\exp_y\bigl(\mathcal B_y(r_x)\bigr)\subset B(y,1).$

The family $\{U_x:x\in M\}$ is an open cover of $M$. Since $M$ is second countable, there exists
a countable subcover
\[
M=\bigcup_{i=1}^{\infty}U_i,
\]
where each $U_i$ is associated with a radius $r_i\in(0,1)$ having
all the properties established above. Define a pairwise disjoint Borel partition of $M$ by
\[
A_1=U_1,
\qquad
A_i=U_i\setminus\bigcup_{j=1}^{i-1}U_j,
\quad i\geq2,
\]
and define the  regular radius as $\rho_b(y)=r_i$ whenever $y\in A_i$. Since each $A_i$ is Borel, the function $\rho_b$ is Borel
measurable. If $y\in A_i$, then $y\in U_i$, so all the required
properties hold on $\mathcal B_y(\rho_b(y))
=\mathcal B_y(r_i)$. This proves the lemma.
\end{proof}

We next adapt the regular radius to a $C^2$ endomorphism $f$ and estimate the differential of its exponential coordinate representation.
\begin{lemma}
\label{26-7-26-1607}
Fixing  $b\in(0, 1/100)$, let $\rho_b: M\to (0,1)$ be the regular radius defined as \Cref{26-7-26-1608}.
Then, for every $C^2$ endomorphism $f:M\to M$ and every $x\in M$,  one has that
\[
f\bigl(
\exp_x(\mathcal B_x(\rho_{b}(x)L(x)^{-1})\rho_b(f(x)))
\bigr)
\subset
\exp_{f(x)}
\bigl(
\mathcal B_{f(x)}(\rho_b(f(x)))
\bigr),
\]
where $L(x):=\max\{4\sup_{y\in B(x,1)}|\mathrm{d}_yf|, 1\}$.  Moreover, for every
$v,w\in\mathcal B_x(\rho_{b}(x)L(x)^{-1})\rho_b(f(x))$, one has
\begin{align*}
&
\left|
\mathrm d_v
\bigl(
\exp_{f(x)}^{-1}\circ f\circ\exp_x
\bigr)
-
\mathrm d_w
\bigl(
\exp_{f(x)}^{-1}\circ f\circ\exp_x
\bigr)
\right|
\\
&\quad\leq
2\Big(
\sup_{y\in B(x,1)}|\mathrm d_yf|
+
\sup_{y\in B(x,1)}|\nabla_y\mathrm df|
+\big(
\sup_{y\in B(x,1)}|\mathrm d_yf|
\big)^2
\Big)|v-w|.
\end{align*}

\end{lemma}

\begin{proof}

Since $M$ is complete, the Hopf--Rinow theorem implies that
$\overline{B(x,1)}$ is compact for any $x\in M$.  For convenience,  denote 
$$\rho_{b,f}(x):=\rho_b(x)L(x)^{-1}\rho_b(f(x)).$$
Let $v\in\mathcal B_x(\rho_{b,f}(x))$
and consider the radial geodesic
\[
\gamma(t)=\exp_x(tv),\qquad 0\leq t\leq1.
\]
Since $|v|<\rho_{b,f}(x)\leq\rho_b(x)<1$,
we have $\gamma(t)\in B(x,1)$ for every $t\in[0,1]$. Moreover,
$\gamma$ has constant speed $|v|$. Hence
\begin{align*}
d\bigl(f(\exp_x(v)),f(x)\bigr)
&\leq
\int_0^1
\left|
\mathrm d_{\gamma(t)}f
\bigl(\dot\gamma(t)\bigr)
\right|\,\mathrm dt
\\
&\leq \sup_{y\in B(x,1)}|\mathrm{d}_yf|\cdot|v|
\\
&< \sup_{y\in B(x,1)}|\mathrm{d}_yf|\cdot\rho_b(x)L(x)^{-1}\rho_b(f(x))
\\&\leq \frac{1}{4}\rho_b(f(x)).
\end{align*}
By completeness, every point of
$B(f(x),\rho_b(f(x)))$ can be joined to $f(x)$ by a minimizing
geodesic. Therefore, $B(f(x),\rho_b(f(x)))\subset\exp_{f(x)}\bigl(\mathcal B_{f(x)}(\rho_b(f(x)))\bigr),$
which implies that
\[
f\bigl(
\exp_x(\mathcal B_x(\rho_{b,f}(x)))
\bigr)
\subset
\exp_{f(x)}
\bigl(
\mathcal B_{f(x)}(\rho_b(f(x)))
\bigr).
\]
Thus the map
\[
F_x:=
\exp_{f(x)}^{-1}\circ f\circ\exp_x
\]
is well defined on $\mathcal B_x(\rho_{b,f}(x))$.

 Fix $u\in\mathcal B_x(\rho_{b,f}(x))$ and set $z=\exp_x(u).$
Since $\rho_{b,f}(x)\leq\rho_b(x)<1$, we have $z\in B(x,1)$.
The identity $\exp_{f(x)}\circ F_x=f\circ\exp_x$
implies
\[
\mathrm d_{F_x(u)}\exp_{f(x)}
\circ\mathrm d_uF_x
=
\mathrm d_zf\circ\mathrm d_u\exp_x.
\]
By the assumptions on the exponential maps, one has that 
\[
\left|
\bigl(
\mathrm d_{F_x(u)}\exp_{f(x)}
\bigr)^{-1}
\right|
\leq\frac{1}{1-b}
\]
and $|\mathrm d_u\exp_x|\leq1+b.$
It follows that
\begin{align}
\label{26-7-26-1721}
|\mathrm d_uF_x|\leq\frac{1+b}{1-b}\sup_{y\in B(x,1)}|\mathrm d_yf|.
\end{align}
Taking the covariant second derivative of $\exp_{f(x)}\circ F_x=f\circ\exp_x,$
we obtain, for every $\xi,\eta\in T_xM$,
\begin{align}
&
\label{26-7-26-1722}\mathrm d_{F_x(u)}\exp_{f(x)}
\bigl(
\mathrm \nabla_u \mathrm{d}F_x(\xi,\eta)
\bigr)
\\\notag&\quad=
(\nabla_z\mathrm df)
\bigl(
\mathrm d_u\exp_x(\xi),
\mathrm d_u\exp_x(\eta)
\bigr)+
\mathrm d_zf
\bigl(
(\nabla\mathrm d\exp_x)_u(\xi,\eta)
\bigr)
\\\notag&\qquad\quad-(\nabla\mathrm d\exp_{f(x)})_{F_x(u)}
\bigl(
\mathrm d_uF_x(\xi),
\mathrm d_uF_x(\eta)
\bigr),
\end{align}
where $\nabla_u\mathrm{d}F_x$ denotes the ordinary second derivative of the
map between the vector spaces $T_xM$ and $T_{f(x)}M$. According to \eqref{26-7-26-1721} and \eqref{26-7-26-1722}, we obtain

\begin{align}
\label{26-7-26-1727}|\nabla_u\mathrm{d}F_x|
\leq{}&
\frac{(1+b)^2}{1-b}
\sup_{y\in B(x,1)}|\nabla_y\mathrm df|+
\frac{b}{1-b}
\sup_{y\in B(x,1)}|\mathrm d_yf|
\\\notag&+
\frac{b(1+b)^2}{(1-b)^3}
\left(
\sup_{y\in B(x,1)}|\mathrm d_yf|
\right)^2.
\end{align}
Since $0<b\leq1/100$, one has that 
\begin{align*}
|\nabla_u\mathrm{d}F_x|\leq2\Big(&\sup_{y\in B(x,1)}|\mathrm d_yf|+\sup_{y\in B(x,1)}|\nabla_y\mathrm df|+\big(\sup_{y\in B(x,1)}|\mathrm d_yf|\big)^2\Big).
\end{align*}
Finally, for every $v,w\in\mathcal B_x(\rho_{b,f}(x))$ one has 
\[
\mathrm d_vF_x-\mathrm d_wF_x=\int_0^1\nabla_{w+t(v-w)}\mathrm{d}F_x
(v-w,\cdot)\,\mathrm dt.
\]
It follows that 
\begin{align*}
|\mathrm d_vF_x-\mathrm d_wF_x|
\leq
2\Big(
&
\sup_{y\in B(x,1)}|\mathrm d_yf|
+
\sup_{y\in B(x,1)}|\nabla_y\mathrm df|+
\big(
\sup_{y\in B(x,1)}|\mathrm d_yf|
\big)^2
\Big)|v-w|.
\end{align*}
This completes the proof.
\end{proof}

 \subsection*{Acknowledgments and Statements}
Zeng Lian was supported by  the NSFC of China  (Grant  No. 12090010) and National Tianyuan Mathematics Southwest Center (Grant  No. 12226102).  Jianhua Zhang was supported by  the China Postdoctoral Science Foundation (Grant  No. 2025M774190). This work is original research and has not been submitted to any other journals for publication. We declare that we do not have any commercial or associative interest that represents a conflict of interest in connection with the work submitted.

\providecommand{\bysame}{\leavevmode\hbox to3em{\hrulefill}\thinspace}
\providecommand{\MR}{\relax\ifhmode\unskip\space\fi MR }
\providecommand{\MRhref}[2]{%
  \href{http://www.ams.org/mathscinet-getitem?mr=#1}{#2}
}
\providecommand{\href}[2]{#2}

\end{document}